# Fourth-Order Paired-Explicit Runge-Kutta Methods

Daniel Doehring [a,*], Lars Christmann [a,d], Michael Schlottke-Lakemper [a,b,c], Gregor J. Gassner [d], and Manuel Torrilhon [a]

[a] *Applied and Computational Mathematics, RWTH Aachen University, Germany.*
[b] *High-Performance Computing Center Stuttgart (HLRS), University of Stuttgart, Germany.*
[c] *High-Performance Computing, Center for Advanced Analytics and Predictive Sciences, University of Augsburg, Germany.*
[d] *Department of Mathematics and Computer Science, Center for Data and Simulation Science, University of Cologne, Germany.*

**Abstract**

In this paper, we extend the Paired-Explicit Runge-Kutta schemes by Vermeire et. al. [1, 2] to fourth-order of consistency. Based on the order conditions for partitioned Runge-Kutta methods we motivate a specific form of the Butcher arrays which leads to a family of fourth-order accurate methods. The employed form of the Butcher arrays results in a special structure of the stability polynomials, which needs to be adhered to for an efficient optimization of the domain of absolute stability.

We demonstrate that the constructed fourth-order Paired-Explicit Runge-Kutta methods satisfy linear stability, internal consistency, designed order of convergence, and conservation of linear invariants. At the same time, these schemes are seamlessly coupled for codes employing a method-of-lines approach, in particular without any modifications of the spatial discretization. We apply the multirate Paired-Explicit Runge-Kutta (P-ERK) schemes to inviscid and viscous problems with locally varying wave speeds, which may be induced by non-uniform grids or multiscale properties of the governing partial differential equation. Compared to state-of-the-art optimized standalone methods, the multirate P-ERK schemes allow significant reductions in right-hand-side evaluations and wall-clock time, ranging from 40% up to factors greater than three.

## 1. Introduction

The numerical solution of time-dependent, convection-dominated nonlinear partial differential equations (PDEs) describing flow phenomena is typically done via explicit time stepping. This introduces a timestep restriction due to the Courant-Friedrichs-Lewy (CFL) condition which is based on the fastest wave speed present in the domain. In the case of locally varying wave speeds due to non-uniform grids or multiscale phenomena of the governing PDE global time-stepping methods are not optimally efficient. Thus, multirate [1–7] and local-time stepping methods [8–14] have been developed to allow for effectively different timesteps in different regions of the computational domain. This goes back to the 1980s when Osher and Sanders devised one of the first multirate methods [6], which has been extended to higher order of accuracy [4, 3, 7, 13]. Around the same time, Berger and Oliger introduced local time stepping in combination with adaptive mesh refinement on Cartesian grids [8]. Since then, many related methods have been developed, including specific techniques for the Discontinuous Galerkin (DG) method [11, 12, 10] or specific problems such as wave-propagation [9]. Recently, this topic received contributions from the field of multirate time integration [1, 2] and classic local time stepping [14].

[*] Corresponding author

While implicit time-stepping methods relieve the timestep restriction of explicit methods, they introduce the necessity of solving a nonlinear system of equations at each timestep or even at each Runge-Kutta stage. Standard Newton-based nonlinear solvers require the Jacobian of the right-hand-side (RHS), which is very memory intensive if stored in dense format for large-scale simulations. While this issue may be addressed by encoding the sparsity of the system, parallelization of the involved linear solver poses another challenge. Additionally, the issue of large truncation errors for large timesteps (which are required to render implicit methods efficient) might arise.

In this work, we consider hyperbolic-parabolic PDEs of the form

$$\partial_t \boldsymbol{u}(t,\boldsymbol{x}) + \nabla \cdot \boldsymbol{f}(\boldsymbol{u}(t,\boldsymbol{x}), \nabla \boldsymbol{u}(t,\boldsymbol{x})) = \boldsymbol{s}(t,\boldsymbol{x},\boldsymbol{u}(t,\boldsymbol{x})) \tag{1.1}$$

with flux function $\boldsymbol{f}$, which may depend on viscous gradients $\nabla \boldsymbol{u}$. Given the inherently different nature of time (directed, one dimensional) and space (non-directed, multiple dimensions) the temporal and spatial derivatives in (1.1) are oftentimes discretized separately. This approach is commonly referred to as the method of lines (MoL) [15], where the application of spatial-derivative discretization techniques such as Finite Differences/Volumes/Elements or DG leads to a potentially huge system of ordinary differential equations (ODEs). We consider the resulting initial value problem (IVP)

$$\boldsymbol{U}(t_0) = \boldsymbol{U}_0 \tag{1.2a}$$

$$\boldsymbol{U}'(t) = \boldsymbol{F}(t, \boldsymbol{U}(t)). \tag{1.2b}$$

with $\boldsymbol{U}, \boldsymbol{F} \in \mathbb{R}^M$. A benefit of the MoL approach is that the PDE can be solved in a modular fashion, i.e., the spatial discretization can be performed independently of the time integration and vice versa. Besides the ease of implementation, the MoL approach benefits from the availability of a large body of literature on both the spatial discretization techniques and time integration methods. This enables an independent analysis of the building blocks of the fully discrete system, see for instance the review paper [15] for the Runge-Kutta Discontinuous Galerkin (RKDG) method.

Recently, a class of stabilized optimized Partitioned Runge-Kutta methods (PRKMs) have been proposed [1, 2] which achieve multirate effects while being readily implemented, as no changes in the spatial discretization are required. Termed Paired-Explicit Runge-Kutta (PERK) [1], these schemes amount only to applying different Runge-Kutta methods to the components of the semidiscretization (1.2) with no need for specialized inter-/extrapolation of interface values or sophisticated flux updating strategies. In fact, they can be implemented similar to the well-known implicit-explicit (IMEX) schemes [16], which apply different Runge-Kutta methods for the different operators present in the numerical scheme. Furthermore, as the P-ERK schemes can be cast into the framework of PRKMs, one can utilize available results on order conditions [17], conservation of linear invariants [18], and nonlinear stability properties [19–22].

In this work, we extend the P-ERK schemes to fourth-order accuracy. We begin in Section 2 with a brief review of PRKMs and conditions for convergence order, absolute stability, conservation, internal consistency, and nonlinear stability. Then, we introduce the P-ERK schemes in Section 3 and extend them to fourth-order in Section 3.2. Due to the form of the order conditions for fourth-order PRKMs the optimization of the domain of absolute stability of the Runge-Kutta methods needs to be revised, which is addressed in Section 4. In Section 5 we demonstrate that the constructed fourth-order P-ERK satisfy linear stability, internal consistency, order of convergence, and conservation of linear invariants. Applications to a wide range of inviscid and viscous flow problems are presented in Section 6 and we conclude the paper in Section 7.

## 2. Partitioned Runge-Kutta Methods

We consider coefficient-based partitioned ODEs of the form

$$\boldsymbol{U}(t_0) = \boldsymbol{U}_0 \tag{2.1a}$$

$$\boldsymbol{U}'(t) = \begin{pmatrix} \boldsymbol{U}^{(1)}(t) \\ \vdots \\ \boldsymbol{U}^{(R)}(t) \end{pmatrix}' = \begin{pmatrix} \boldsymbol{F}^{(1)}(t, \boldsymbol{U}^{(1)}(t), \dots, \boldsymbol{U}^{(R)}(t)) \\ \vdots \\ \boldsymbol{F}^{(R)}(t, \boldsymbol{U}^{(1)}(t), \dots, \boldsymbol{U}^{(R)}(t)) \end{pmatrix} = \boldsymbol{F}(t, \boldsymbol{U}(t)) \tag{2.1b}$$

where $R \in \mathbb{N}$ denotes the number of partitions/levels. In the case of non-uniform meshes, the partitioning of the semidiscretization (1.2) is based on the minimum edge length $h$ of the grid cells. In the context of the DG method employed in this work, the coefficients of the local solution polynomials are the unknowns of the ODE system (2.2) and can be uniquely assigned to a partition based on the minimum edge length $h$ of the associated cell. For systems with strongly spatially varying eigenvalues of the flux-Jacobians $J_i := \partial_{\boldsymbol{u}} \boldsymbol{f}_i, i = 1, 2, 3$, the partitioning of the ODE system (2.2) may be based on the magnitude of the spectral radii $\rho_i$ of the flux Jacobians $J_i$. We present applications of the P-ERK methods for both of these scenarios.

In the following, the superscript $(\cdot)^{(r)}$ indicates quantities corresponding to the $r$'th partition. These may be a subset of unknowns $\boldsymbol{U}^{(r)} \in \mathbb{R}^{M_r}$, or a method used for solving the corresponding system with RHS $\boldsymbol{F}^{(r)}$. ODE systems of the form (2.2) encourage the use of PRKMs, which, in Butcher form, are given by [23, Chapter II.15]

$$\boldsymbol{U}_0 = \boldsymbol{U}(t_0)\,, \tag{2.2a}$$

$$\boldsymbol{K}_i^{(r)} = \boldsymbol{F}^{(r)}\left(t_n + c_i^{(r)}\Delta t, \boldsymbol{U}_n + \Delta t \sum_{j=1}^{S}\sum_{k=1}^{R} a_{i,j}^{(k)} \boldsymbol{K}_j^{(k)}\right), \quad i = 1, \dots, S;\;\; r = 1, \dots, R \tag{2.2b}$$

$$\boldsymbol{U}_{n+1} = \boldsymbol{U}_n + \Delta t \sum_{i=1}^{S}\sum_{r=1}^{R} b_i^{(r)} \boldsymbol{K}_i^{(r)}\,, \tag{2.2c}$$

where $S$ denotes the number of stages $\boldsymbol{K}_i \in \mathbb{R}^M$.

PRKMs have been introduced already in the 1970s [24, 25] and have received notable usage for the integration of Hamiltonian systems, see [26, 27] and references therein. An essentially equivalent class of Runge-Kutta methods are Additive Runge-Kutta methods (ARKMs) [28–30] which have been developed to achieve efficient integration of multiscale and stiff systems [5, 29]. Similar works have been performed with PRKMs [3, 31, 32]. Since the P-ERK methods we apply in this work are a special case of PRKMs [2], we now discuss conditions for order of consistency, absolute stability, conservation of linear invariants, internal consistency, and nonlinear stability.

### *2.1. Order Conditions*

On top of the classical conditions for convergence of order $p$, PRKMs need to satisfy additional conditions which have been derived by multiple authors, see [25, 17, 33] and for a more recent, compact summary [34]. For identical abscissae $\boldsymbol{c}^{(r)} \equiv \boldsymbol{c}$ and weights $\boldsymbol{b}^{(r)} \equiv \boldsymbol{b}$ (more on this restriction in sections 2.3, 2.4), these in principle additional order constraints reduce for $p = 1, 2, 3$ conveniently to the classical ones [34], [35, Chapter 5.9]:

$$p = 1: \qquad \boldsymbol{b}^T \boldsymbol{1} \overset{!}{=} 1 \tag{2.3a}$$

$$p = 2: \qquad \boldsymbol{b}^T \boldsymbol{c} \overset{!}{=} \frac{1}{2} \tag{2.3b}$$

$$p = 3: \qquad \boldsymbol{b}^T \boldsymbol{c}^2 \overset{!}{=} \frac{1}{3} \tag{2.3c}$$

$$\boldsymbol{b}^T A^{(r)} \boldsymbol{c} \overset{!}{=} \frac{1}{6} \quad \forall\, r = 1, \dots, R. \tag{2.3d}$$

Here, $\mathbf{1} \in \mathbb{R}^S$ denotes the column-vector of ones and the exponentiation of vectors is to be understood element-wise. For order $p = 4$, however, even for identical abscissae $\boldsymbol{c}$ and weights $\boldsymbol{b}$, non-trivial coupling conditions between the Butcher arrays $A^{(r)}$ of the different methods arise [34]:

$$\boldsymbol{b}^T \boldsymbol{c}^3 \overset{!}{=} \frac{1}{4} \tag{2.4a}$$

$$\boldsymbol{b}^T C A^{(r)} \boldsymbol{c} \overset{!}{=} \frac{1}{8} \quad \forall\, r = 1, \dots, R \tag{2.4b}$$

$$\boldsymbol{b}^T A^{(r)} \boldsymbol{c}^2 \overset{!}{=} \frac{1}{12} \quad \forall\, r = 1, \dots, R \tag{2.4c}$$

$$\boldsymbol{b}^T A^{(r_1)} A^{(r_2)} \boldsymbol{c} \overset{!}{=} \frac{1}{24} \quad \forall\, r_1, r_2 = 1, \dots, R \tag{2.4d}$$

rendering the construction of an optimized fourth-order P-ERK scheme significantly more difficult. Here, $C \in \mathbb{R}^{S\times S}$ denotes the diagonal matrix constructed from the abscissae, i.e., $C := \operatorname{diag}(\boldsymbol{c})$.

While (2.4a)–(2.4d) reduce for $r_1 = r_2$ to the usual fourth-order conditions for the $r$'th method [35, Chapter 5.9], i.e., impose only constraints on the $r$'th Butcher array $A^{(r)}$, the coupling conditions (2.4d) demand a compatibility of the Butcher arrays $A^{(r_1)}$ and $A^{(r_2)}$ for $r_1 \neq r_2$. This leaves one with two options: Either optimize the entire family of Butcher arrays $A^{(r)}, r = 1, \dots, R$ simultaneously or make a choice for the functional form of the Butcher arrays $A^{(r)}$ which satisfies (2.4d) by construction and optimize the remaining free parameters. While the first option offers in principle larger optimal timesteps, it comes with the drawback of a significantly more complex optimization problem, as there are many more free parameters which need to be optimized jointly. Additionally, the extension of such a family by another scheme is in general not possible. In contrast, by choosing a specific form of the Butcher arrays $A^{(r)}$ which satisfies (2.4d) by construction, one can optimize each method separately, thus drastically simplifying the optimization problem and allowing for easy extension of the family by another scheme. Consequently, a special form of the Butcher arrays $A^{(r)}$ is sought satisfying (2.4d) while still allowing optimization of the associated stability polynomial $P^{(r)}(z)$. We propose such an archetype method in Section 3.2.

### *2.2. Absolute Stability*

Absolute/linear stability examines the asymptotic behaviour of a time integration method applied to the constant-coefficient linear test system [36, Chapter IV.2]

$$\boldsymbol{U}'(t) = J\,\boldsymbol{U}(t)\,, \tag{2.5}$$

where $J$ denotes the Jacobian of the RHS $\boldsymbol{F}(t, \boldsymbol{U}(t))$ (1.2)

$$J(t, \boldsymbol{U}) := \frac{\partial \boldsymbol{F}(t, \boldsymbol{U}(t))}{\partial \boldsymbol{U}}\,. \tag{2.6}$$

A time integration method is called linearly/absolutely stable if its associated stability function $P(z)$ (which is a polynomial for explicit methods) is in magnitude less or equal than one when evaluated at the damping modes, i.e., scaled eigenvalues with non-positive real part

$$|P(\Delta t \lambda)| \overset{!}{\leq} 1 \quad \forall\, \lambda \in \boldsymbol{\sigma}(J) : \operatorname{Re}(\lambda) \leq 0\,. \tag{2.7}$$

In other words, (2.7) demands that the scaled spectrum lies in the region of absolute stability of the method with stability polynomial $P(z)$. Since the RHS $\boldsymbol{F}$ depends on the solution $\boldsymbol{U}$ itself, we require (2.7) to hold for all (unknown) states $\boldsymbol{U}(t)$ reached over the course of a simulation. This is in principle difficult to assure, but the spectra $\boldsymbol{\sigma}(J(\boldsymbol{U}))$ are observed to be in practice relatively robust, i.e., similar to the spectrum of the initial condition $\boldsymbol{\sigma}\left(J(\boldsymbol{U}_0)\right)$. With 'similar' we refer to the characteristic distribution of the eigenvalues in the complex plane, which remains robust over the course of a simulation (up to scaling due to mesh refinement), even for discontinuous solutions. In fact, we have observed that the spectrum is mainly influenced by the simulated PDE, number of spatial dimensions, local solution polynomial degree $k$, and approximate Riemann solver. In comparison, initial and boundary conditions have only a minor impact on the spectrum. Exemptions from this are for instance simulations where the medium is initially at rest and e.g. a supersonic jet is prescribed at a boundary.

The linear stability properties of PRKMs have been studied in [37] with special focus on separable Hamiltonian systems. For the related ARKMs, a linear stability analysis is performed in [30] based on the test problem $U' = \sum_{r=1}^{R} \lambda^{(r)} U$.
For PRKMs it is customary to investigate the test equation

$$\boldsymbol{U}' = \Lambda \boldsymbol{U}, \quad \Lambda \in \mathbb{C}^{N \times N}, \tag{2.8}$$

which can be naturally partitioned according to a set of mask matrices $I^{(r)} \in \{0,1\}^{N \times N}, r = 1, \dots, R$:

$$\Lambda = \sum_{r=1}^{R} \Lambda^{(r)} = \sum_{r=1}^{R} I^{(r)} \Lambda . \tag{2.9}$$

The mask matrices sum up to the identity matrix $I$ of dimensions $N \times N$, i.e., $\sum_{r=1}^{R} I^{(r)} = I$ and are required to specify the exact partitioning of (2.1a):

$$\boldsymbol{U}(t) = \sum_{r=1}^{R} \underbrace{I^{(r)} \boldsymbol{U}(t)}_{=:\boldsymbol{U}^{(r)}(t)}, \qquad \boldsymbol{F}(t, \boldsymbol{U}(t)) = \sum_{r=1}^{R} \underbrace{I^{(r)} \boldsymbol{F}(t, \boldsymbol{U}(t))}_{=:\boldsymbol{F}^{(r)}\left(t, \boldsymbol{U}(t)\right)} . \tag{2.10}$$

The matrix-valued stability function $P(Z) \in \mathbb{C}^{N \times N}$ with scaled system matrix

$$Z^{(r)} \coloneqq \Delta t \Lambda^{(r)} \tag{2.11}$$

is given by [18]

$$P(Z) = I + \left( \sum_{r=1}^{R} \left( \boldsymbol{b}^{(r)} \otimes I \right)^T \left( Z^{(r)} \otimes I \right) \right) \left( I_{NS} - \sum_{r=1}^{R} \left( A^{(r)} \otimes I \right) \left( Z^{(r)} \otimes I \right) \right)^{-1} (\boldsymbol{1} \otimes I) \tag{2.12}$$

where $\otimes$ denotes the Kronecker product and $I_{NS}$ the identity matrix of dimension $N \cdot S$. As discussed in [37], it is difficult to infer from (2.12) information on the linear stability of the overall partitioned scheme based on the stability properties of the individual schemes $(A^{(r)}, b^{(r)})$ for arbitrary mask matrices $I^{(r)}$. In contrast to classic Runge-Kutta methods for PRKMs one can no longer solely infer the action of the time integration scheme from the stability functions of the individual methods [37].

In practice, however, it is observed that absolute stability of each method $(A^{(r)}, \boldsymbol{b}^{(r)})$ for its associated spectrum $\boldsymbol{\lambda}^{(r)}$ suffices for overall linear stability, irrespective of the concrete realization of the mask matrices $I^{(r)}$ [1, 2]. To check this, one can either compute $P(Z)$ via (2.12) or set up the fully discrete system directly, i.e., apply the Runge-Kutta method to (2.5). Then, linear stability follows from applying stability criteria for discrete-time linear time-invariant systems [38, Chapter 8.6]. The latter approach is illustrated for a concrete example in Section 5.1.

### *2.3. Internal Consistency*

To ensure that the approximations of the partitioned stages $\boldsymbol{K}_i^{(r)}$ approximate indeed the same timesteps $c_i$ across partitions $r$ we require internal consistency [18, 21, 30], i.e.,

$$\sum_{j=1}^{i-1} a_{i,j}^{(r)} \overset{!}{=} c_i \quad \forall\, i = 1, \dots, S; \ \forall\, r = 1, \dots, R \tag{2.13}$$

which is equivalent to all $R$ Runge-Kutta methods having stage order one [18]. Thus, $\boldsymbol{c}^r \equiv \boldsymbol{c},\ \forall\, r = 1, \dots, R$. It has been demonstrated in [21, 18] that internal consistency is required to avoid spurious oscillations at the interfaces between partitions. We demonstrate the internal consistency of the P-ERK schemes by comparing the errors of a family of methods ($R > 1$) to each composing individual scheme in Section 5.2.

### *2.4. Conservation of Linear Invariants*

Since the underlying PDEs, from which the semidiscretization (1.2) is constructed, correspond oftentimes to conservation laws (i.e., (1.1) with no sources $\boldsymbol{s} \equiv \boldsymbol{0}$) prescribing the conservation of mass, momentum, and energy, it is natural to require that the application of a PRKM preserves this property for conservative spatial discretizations such as Finite Volume and DG schemes. For equation-based partitioning according to (2.1a) this is ensured if the weights $\boldsymbol{b}$ are identical across partitions [4, 18], i.e.,

$$b_i^{(r_1)} = b_i^{(r_2)} = b_i, \quad \forall\, i = 1, \dots, S;\ \forall\, r_1, r_2 = 1, \dots, R\,. \tag{2.14}$$

For flux-based partitioning this restriction may be relaxed, see [22] for a detailed discussion. An experimental confirmation of the conservation of linear invariants is provided in Section 5.3.

### *2.5. Nonlinear Stability*

Nonlinear stability includes (among stronger results for scalar one-dimensional equations such as the total-variation diminishing (TVD) property) positivity preservation of physical quantities such as pressure and density and the suppression of spurious oscillations around discontinuities. Time integration methods guaranteeing these properties are coined strong stability preserving (SSP) [39] and have been well-studied over the past years, see [40] and references therein. The nonlinear stability properties of partitioned and additive Runge-Kutta methods have been thoroughly investigated in [19–22]. In the aforementioned works it was shown that for an overall SSP PRKM the individual methods have to be SSP and the nonlinearly stable timestep $\Delta t$ is restricted by

$$\Delta t = \min_r \left\{ c_r^{\text{SSP}} \right\} \Delta t_{\text{Forward Euler}}. \tag{2.15}$$

Here, $c_r^{\text{SSP}}$ denotes the SSP coefficient [40, Chapter 3] (radius of absolute monotonicity [41, 19]) of the $r$'th method. Except for some special cases where the timestep is restricted due to absolute stability (see e.g. [42]), (2.15) renders the application of SSP-capable PRKMs in many cases ineffective, as the method with the smallest SSP coefficient governs the overall admissible timestep which guarantees nonlinear stability. For the P-ERK schemes, we conducted an in-depth nonlinear stability analysis in [43] and provide a related discussion involving a numerical example in Appendix B.

Having discussed the properties of PRKMs, we now turn our attention to a special case thereof, namely the Paired-Explicit Runge-Kutta (P-ERK) methods proposed by Vermeire [1].

## 3. Paired-Explicit Runge-Kutta (P-ERK) Methods

Paired-Explicit Runge-Kutta (P-ERK) methods have been originally proposed in [1] and have been extended to third-order in [2] and embedded schemes with error estimation in [44]. In the aforementioned works, the P-ERK schemes have been applied to solve compressible Navier-Stokes equations on non-uniform grid, where speedup is gained due to reduced computational effort on the coarse grid cells. In a previous work, we applied the P-ERK schemes to hyperbolic-parabolic and purely hyperbolic problems on adaptively refined meshes [43].

### *3.1. Central Idea*

It is well-known that explicit methods are only conditionally stable for IVPs (1.2) with time step restriction due to some form of Courant-Friedrichs-Lewy (CFL) condition:

$$\Delta t \overset{!}{\leq} \text{CFL} \cdot \frac{1}{\nu}. \tag{3.1}$$

Here $\nu > 0$ denotes an (estimate) for the largest characteristic speed. For the DG method employed in this work, the local (inviscid) wave speed $\nu$ scales as

$$\nu \sim \max_{i=1,\dots,N_D} \frac{(k+1) \cdot \rho_i}{h_i} \tag{3.2}$$

where $N_D$ denotes the number of spatial dimensions, $h_i$ is the smallest characteristic cell size in $i$-direction, $k$ is the solution polynomial degree, and $\rho_i$ is the spectral radius of the directional flux Jacobian

$$J_i := \frac{\partial \boldsymbol{f}_i}{\partial \boldsymbol{u}}. \tag{3.3}$$

For non-uniform grids, i.e., varying $h_i$, the CFL condition (3.1) restricts the time step size to the smallest cell size $h_i$. When using a traditional explicit method, this causes unnecessary computational effort for the larger cells. To overcome this issue, the P-ERK methods employ Runge-Kutta schemes with higher CFL numbers in regions with high $\nu$ and schemes with lower CFL numbers in regions with low $\nu$. In particular, optimized Runge-Kutta methods with different domains of absolute stability are combined to form a P-ERK family.

The estimate (3.2) is only a valid bound for the maximum timestep wave speed in the convection-dominated case, i.e., when the spectral radius $\rho_{\boldsymbol{F}}$ of the Jacobian of the fully discrete system (1.2b) $J_{\boldsymbol{F}}(\boldsymbol{U}) := \frac{\partial \boldsymbol{F}}{\partial \boldsymbol{U}}$ scales as

$$\rho \sim \frac{|a|}{h} \gg \frac{|d|}{h^2}\,, \tag{3.4}$$

where $h$ denotes the smallest characteristic cell size and $a, d$ are suitable quantifiers for the influence of convection and diffusion, respectively.

For diffusion-dominated PDEs, the spectral radius increases quadratically with the inverse of the minimal grid size $h$, cf. (3.4). Thus, for an efficient treatment of such problems, the region of stability needs to increase also quadratically with the number of stages evaluations $E$, which is the case for Runge-Kutta-Chebyshev methods [45, 46] targeting precisely parabolic problems. These methods, however, do not extend (except for special treatment [47, 48]) significantly into the complex plane, rendering them unstable for the convection-induced complex eigenvalues.

In the convection-dominated case, however, the characteristic wave speed scales linearly inversely proportional to the smallest cell size $h$, cf. (3.2). This is crucial for the effectiveness of the P-ERK schemes, as the maximum stable timestep of optimized Runge-Kutta methods scales asymptotically linearly with the number of stage evaluations $E$ for hyperbolic spectra [49–51]. Consequently, the usage of an optimized $2E$ stage-evaluation Runge-Kutta method for the simulation of a convection-dominated PDE with smallest grid size $\frac{1}{2}h$ allows using the same timestep $\Delta t$ as the simulation of the same problem on a grid with smallest size $h$ with an $E$ stage-evaluation scheme.

For illustration of the idea we present the Butcher tableaus of an second-order, $S = 6$ stage,

$E = \{3, 6\}$ stage evaluation P-ERK method:

$$
\begin{array}{c|c|cccccc|cccccc}
i & \boldsymbol{c} & & & A^{(1)} & & & & & & A^{(2)} & & & \\
\hline
1 & 0 & & & & & & & & & & & & \\
2 & 1/10 & 1/10 & & & & & & 1/10 & & & & & \\
3 & 2/10 & 2/10 & 0 & & & & & 2/10 - a^{(2)}_{3,2} & a^{(2)}_{3,2} & & & & \\
4 & 3/10 & 3/10 & 0 & 0 & & & & 3/10 - a^{(2)}_{4,3} & 0 & a^{(2)}_{4,3} & & & \\
5 & 4/10 & 4/10 & 0 & 0 & 0 & & & 4/10 - a^{(2)}_{5,4} & 0 & 0 & a^{(2)}_{5,4} & & \\
6 & 5/10 & 5/10 - a^{(1)}_{6,5} & 0 & 0 & 0 & a^{(1)}_{6,5} & & 5/10 - a^{(2)}_{6,4} & 0 & 0 & 0 & a^{(2)}_{6,4} & \\
\hline
 & \boldsymbol{b}^T & 0 & 0 & 0 & 0 & 0 & 1 & 0 & 0 & 0 & 0 & 0 & 1
\end{array}
\tag{3.5}
$$

The first method $A^{(1)}$ requires only the computation of stages $\boldsymbol{K}_1^{(1)}, \boldsymbol{K}_5^{(1)}, \boldsymbol{K}_6^{(1)}$, i.e., three evaluations of $\boldsymbol{F}^{(1)}$, while the second method $A^{(2)}$ requires computation of all six stages $\boldsymbol{K}_i^{(2)}, i = 1, \dots, 6$. When looking at the first method in isolation, it becomes a reducible [36] method, i.e., the stages $\boldsymbol{K}_2^{(1)}, \boldsymbol{K}_3^{(1)}$ do not influence the final solution $\boldsymbol{U}_{n+1}^{(1)}$ and the Butcher tableau could be truncated to a three-stage method. In the context of PRKMs, however, the second to fourth stage are required for an internally consistent (cf. (2.13)) update of the intermediate state $\boldsymbol{U}_n + \Delta t \sum_{j=1}^S \sum_{r=1}^R a_{i,j}^{(r)} \boldsymbol{K}_j^{(r)}$, see (2.2b).

This particular form of the Butcher tableaus (3.7) comes with three advantages: First, the P-ERK methods are low-storage, requiring only storage of two Runge-Kutta stages $\boldsymbol{K}_i$ at the same time. Second, the computational costs per stage are not increasing if higher stage methods are used, as always at most two previous stages are used for computation of the intermediate state. Third, the sparse structure of both $A^{(r)}$ and $\boldsymbol{b}^{(r)}$ allows for a simplified computation of the coefficients $a_{i,j}^{(r)}$ from an optimized stability polynomial.

Furthermore, the P-ERK schemes offer the flexibility to adapt the coefficients $a_{i,j}^{(r)}$ to the specific problem at hand by optimizing the absolute domain of stability, i.e., maximizing the linearly stable timestep $\Delta t$, cf. (2.7). This optimization is discussed in Section 4 in more detail.

### *3.2. Fourth-Order Paired-Explicit Runge-Kutta Methods*

This section addresses the construction of fourth-order P-ERK methods that fulfill the order conditions (2.3), (2.4) while still allowing for a feasible optimization of the associated stability polynomials $P^{(r)}(z)$ of the individual methods. To start, we need to make some fundamental design choices. As for the third-order P-ERK methods [2] we allow only the last two entries in the weight vector $\boldsymbol{b}$ to be non-zero:

$$\boldsymbol{b}^T = \begin{pmatrix} 0 & \dots & 0 & b_{S-1} & b_S \end{pmatrix}. \tag{3.6}$$

This simplifies the order conditions and sets the stage for a low-storage scheme. Furthermore, we seek to keep the classic P-ERK form of the Butcher array of the different schemes, i.e., having only nonzero entries in the first column and sub-diagonal of $A^{(r)}$:

$$
\begin{array}{c|c|cccccc}
i & \boldsymbol{c} & & & A^{(r)} & & & \\
\hline
1 & 0 & & & & & & \\
2 & c_2 & c_2 & & & & & \\
3 & c_3 & c_3 - a^{(r)}_{3,2} & a^{(r)}_{3,2} & & & & \\
4 & c_4 & c_4 - a^{(r)}_{4,3} & 0 & a^{(r)}_{4,3} & & & \\
\vdots & \vdots & \vdots & & & \ddots & & \\
S & c_S & c_S - a^{(r)}_{S,S-1} & 0 & \dots & 0 & a^{(r)}_{S,S-1} & \\
\hline
 & & 0 & 0 & 0 & 0 & b_{S-1} & b_S
\end{array}
\tag{3.7}
$$

For this particular choice of $\boldsymbol{b}$ and $A^{(r)}$ we can now evaluate the full set of order conditions (2.3), (2.4), which are spelled out in Appendix A. Here, we focus on the coupling condition (2.4d) which for the weight vector (3.6) and Butcher arrays $A^{(r_1)}, A^{(r_2)}$ (3.7) simplifies to

$$\frac{1}{24} \overset{!}{=} b_{S-1}\, a^{(r_1)}_{S-1,S-2}\, a^{(r_2)}_{S-2,S-3}\, c_{S-3} + b_S\, a^{(r_1)}_{S,S-1}\, a^{(r_2)}_{S-1,S-2}\, c_{S-2}, \quad \forall\, r_1, r_2 = 1, \dots, R. \tag{3.8}$$

To satisfy (3.8) for all $r_1, r_2 = 1, \dots, R$ we restrict us to methods which share the same Butcher array coefficients $a_{S-2,S-3}$, $a_{S-1,S-2}$, $a_{S,S-1}$. In addition to the seven variables $a_{S-2,S-3}$, $a_{S-1,S-2}$, $a_{S,S-1}$, $b_{S-1}$, $b_S$, $c_{S-2}$, $c_{S-3}$ from (3.8) we have the timesteps $c_{S-1}$, $c_S$ at our disposal to fulfill the eight order conditions (A.1). So in total there are nine variables and eight order conditions to satisfy, so we expect to find a one-parameter family of fourth-order P-ERK methods. Indeed, using symbolic algebra software we find a set of coefficients that satisfy (A.1) to double precision:

$$b_{S-1} = b_S = 0.5 \tag{3.9a}$$

$$c_{S-2} = 0.479274057836310 \tag{3.9b}$$

$$c_{S-1} = 0.5 + \frac{\sqrt{3}}{6} = 0.788675134594813 \tag{3.9c}$$

$$c_S = 0.5 - \frac{\sqrt{3}}{6} = 0.211324865405187 \tag{3.9d}$$

$$a_{S-2,S-3} = \frac{0.114851811257441}{c_{S-3}} \tag{3.9e}$$

$$a_{S-1,S-2} = 0.648906880894214 \tag{3.9f}$$

$$a_{S,S-1} = 0.0283121635129678 \tag{3.9g}$$

where $c_{S-3}$ is a free parameter. We see from this solution that for a fourth-order accurate P-ERK method we require at least a five stage method as $c_{S-3} \neq 0$ to keep $a_{S-2,S-3}$ finite, cf. (3.9e). This is because we require for an internally consistent method that $\sum_{j=1}^{i-1} a_{i,j} = c_i$, cf. (2.13), which would yield for a four stage method exactly $c_1 = c_{4-3} = 0$. Consequently, by our simplification to consider shared Butcher array coefficients $a_{S-2,S-3}$, $a_{S-1,S-2}$, $a_{S,S-1}$, we pay the price of loosing the possibility to construct a fourth-order P-ERK method with four stages only. This is a qualitative difference to the second and third-order P-ERK methods, which can indeed be constructed with two or three stages, respectively. This is a consequence of our restriction to a two-entry weight vector $\boldsymbol{b}$ enabling a low-storage scheme and the classic P-ERK form of the Butcher array (3.7).

It remains to choose the free parameter $c_{S-3}$ in (3.9). Since for the P-ERK methods we set $a_{i,1} = c_i - a_{i,i-1}$ the stability polynomial $P_{4;5}(z)$ is not affected by changes in $c_{S-3}$, thus other criteria should be employed to choose $c_{S-3}$. Here, we choose $c_{S-3}$ such that the resulting method minimizes the amplification of round-off errors, i.e., optimizes the internal stability properties of the method [52]. Allowing for $c_{S-3} \in [0, 1]$ we find that the optimal value is $c_{S-3} = 1$. A more detailed discussion on internal stability is given in Section 4.4. We emphasize that different criteria for the choice of $c_{S-3}$ are possible, e.g. to minimize dissipation and/or dispersion, see for instance [53, 54] and references therein.

### *3.3. Implementation*

Given that all methods $r = 1, \dots R$ of a family of fourth-order P-ERK schemes share the identical final stages $i = S-3, \dots S$, the fourth-order P-ERK methods are only truly partitioned in the first $i = 3, \dots S-4$ stages. As for the second- and third-order P-ERK schemes, the first and second stage are also shared for the fourth-order P-ERK methods due to the internal consistency requirement (2.13). Thus, the Butcher tableaus $A^{(r)}$ of a range of fourth-order P-ERK methods can be schematically represented as

$$
\begin{array}{c|c|ccc}
 & & & r & \\
\hline
i & \boldsymbol{c} & 1 & \dots & R \\
\cline{3-5}
 & & A^{(1)} & \dots & A^{(R)} \\
\hline
 & & & \boldsymbol{b}^T &
\end{array}
=
\begin{array}{c|c|ccc}
 & & & r & \\
\hline
i & \boldsymbol{c} & 1 & \dots & R \\
\hline
1 & 0 & & & \\
2 & c_2 & & A_{1:2} & \\
\cdashline{3-5}
3 & c_3 & & & \\
\vdots & \vdots & A^{(1)}_{3:S-4} & \dots & A^{(R)}_{3:S-4} \\
S-4 & c_{S-4} & & & \\
\cdashline{3-5}
S-3 & c_{S-3} & & & \\
\vdots & \vdots & & A_{S-3:S} & \\
S & c_S & & & \\
\hline
 & & & \boldsymbol{b}^T &
\end{array}
\tag{3.10}
$$

i.e., they can be seen as a classic-partitioned-classic Runge-Kutta method. Utilizing this leads not only to a more efficient implementation, but also benefits the nonlinear stability properties of the method, see Appendix B.

## 4. Optimization of the Stability Polynomials

In this section we discuss the optimization of the stability polynomials $P^{(r)}_{4;E}(z)$ for the fourth-order P-ERK methods. Due to our choice of the weight vector $\boldsymbol{b}$ (3.6) and the Butcher array $A^{(r)}$ (3.9) the optimization of the stability polynomial becomes more complex compared to the second and third-order case.

### *4.1. Objective of the Optimization*

Recalling the linear stability constraint (2.7) for an $M$-dimensional linear system whose eigenvalues $\lambda_m, m = 1, \dots, M$ all have non-positive real part

$$\left|P_{p;E}(\Delta t \lambda_m)\right| \leq 1, \quad \mathrm{Re}\left(\lambda_m\right) \leq 0, \quad m = 1, \dots, M, \tag{4.1}$$

one is interested in finding the maximum admissible timestep $\Delta t^\star$ and corresponding stability polynomial $P^\star_{p;E}(z)$ of degree $E$ and linear order $p$ corresponding to a certain spectrum $\boldsymbol{\sigma} = \{\lambda_m\}_{m=1,\dots,M}$:

$$\max_{P_{p;E} \in \mathcal{P}_{p;E}} \Delta t \text{ such that } \left|P_{p;E}(\Delta t \lambda_m)\right| \leq 1, \quad m = 1, \dots, M. \tag{4.2}$$

Here, $\mathcal{P}_{p;E}$ denotes the vector space of polynomials of degree $E$ with real coefficients over the complex numbers $\mathbb{C}$ with linear order $p$, i.e., polynomials of form

$$P_{p;E} = \sum_{j=0}^{p} \frac{1}{j!} z^j + \sum_{j=p+1}^{E} \alpha_j z^j, \quad \alpha_j \in \mathbb{R}. \tag{4.3}$$

Given symmetry around the real axis it suffices to check the stability constraint (4.1) only at one eigenvalue of a complex-conjugated pair.

*4.2. Optimization of the Fourth-Order P-ERK Stability Polynomials*

We have seen in Section 3.2 that for the derived family of fourth-order P-ERK methods at least five stages are required. For the five stage-evaluation method the stability polynomial $P_{4;5}(z)$ is fully determined and given by

$$P_{4;5}(z) = \sum_{j=0}^{4} \frac{1}{j!} z^j + k_1 z^5 , \tag{4.4}$$

where the fifth coefficient $k_1 = 0.001055026310046423$ is independent of the choice of $c_{S-3}$. It is instructive to consider the stability polynomial of the eight stage-evaluation method $P_{4;8}(z)$ since it illustrates the relation of the coefficients $\alpha_j$ of the stability polynomial in monomial form, i.e.,

$$P_{4;E}(z; \boldsymbol{\alpha}) = \sum_{j=0}^{4} \frac{1}{j!} z^j + \sum_{j=5}^{E} \alpha_j z^j \tag{4.5}$$

to the free Butcher array coefficients $a_{i,i-1}$. For the sake of readability we define $a_i \coloneqq a_{i,i-1}$. The stability polynomial $P_{4;8}(z)$ is then given by [36, Chapter IV.7]

$$P_{4;8}(z; \boldsymbol{a}) = 1 + z\boldsymbol{b}^T \left(I - zA\right)^{-1} \mathbf{1} \tag{4.6a}$$

$$\begin{aligned} &= \sum_{j=0}^{4} \frac{1}{j!} z^j + \underbrace{\left( \frac{k_2}{c_{8-3}} a_5\, c_4 + k_1 \right)}_{\alpha_5} z^5 + \underbrace{\left( \frac{k_2}{c_{8-3}} a_4\, a_5\, c_3 + \frac{k_1}{c_{8-3}} a_5\, c_4 \right)}_{\alpha_6} z^6 \\ &+ \underbrace{\left( \frac{k_2}{c_{8-3}} a_3\, a_4\, a_5\, c_2 + \frac{k_1}{c_{8-3}} a_4\, a_5\, c_3 \right)}_{\alpha_7} z^7 + \underbrace{\frac{k_1}{c_4} a_3\, a_4\, a_5\, c_2}_{\alpha_8} z^8 \end{aligned} \tag{4.6b}$$

where $k_2 = 0.03726406530405851/c_{S-3}$ is a parameter following from (3.9).

Although maybe not obvious at first sight, there is a clear structure in the coefficients $\alpha_j$ of the stability polynomial $P_{4;8}(z)$ which is crucial for an efficient optimization of the $a_i$. This structure is best seen when comparing the products involving $\frac{k_2}{c_{8-3}}, \frac{k_1}{c_{8-3}}$ of two successive coefficients $\alpha_j, \alpha_{j+1}$. We see from (4.6b) that the factor multiplied by $\frac{k_2}{c_{8-3}}$ in the $\alpha_j$ coefficient becomes exactly the factor multiplied by $\frac{k_1}{c_{8-3}}$ in the $\alpha_{j+1}$ coefficient. Furthermore, the factor multiplied by $\frac{k_2}{c_{8-3}}$ in the $\alpha_{j+1}$ coefficient is the factor multiplied by $\frac{k_2}{c_{8-3}}$ in the $\alpha_j$ coefficient with one additional Butcher array coefficient $a_i$ and exchanged abscissae $c_{i-1}$. Casting this in a general form, the stability polynomials $P_{4;E}(z)$ may be written as

$$\begin{aligned} P_{4;E}(z; \boldsymbol{\gamma}) = &\sum_{j=0}^{4} \frac{1}{j!} z^j + \left( k_1 + \gamma_1^{(E)} \frac{k_2}{c_{E-3}} c_{E-4} \right) z^5 \\ &+ \sum_{k=6}^{E-1} \left( \gamma_{k-4}^{(E)} \frac{k_2}{c_{E-3}} c_{E-k+1} + \gamma_{k-5}^{(E)} \frac{k_1}{c_{E-3}} c_{E-k+2} \right) z^k + \gamma_{E-5}^{(E)} \frac{k_1}{c_{E-3}} c_2 z^E \end{aligned} \tag{4.7}$$

with

$$\gamma_j^{(E)} \coloneqq \prod_{i=1}^{j} a_{E-3-i+1}, \quad j = 1, \dots, E-5 . \tag{4.8}$$

By introducing the $\gamma_k$ we have recovered a linear parametrization of the stability polynomial in the optimization variables $\boldsymbol{\gamma} \in \mathbb{R}_+^{E-5}$, despite using the Butcher-tableau based formulation of the stability polynomial (4.6a). This is a crucial property which allows rewriting the stability constraints (2.7) as convex second-order cone constraints, see [51] and Appendix C. Then, the optimal timestep $\Delta t$ can be found by wrapping the convex optimization problem solver into an bisection routine for the timestep as explained in [51]. After obtaining the optimal $\boldsymbol{\gamma}$ the Butcher array coefficients $a_i$ can be easily computed from (4.8) in a sequential manner.

Note that with the chosen form of the Butcher tableau (3.9), for an $E$-degree stability polynomial there are only $E-5 = E-p-1$ free coefficients $\gamma_j^{(E)}$, as the timesteps $c_i$ are shared among all methods and may not be individually optimized. Again, this is a substantial difference to the second- and third-order P-ERK methods which have $E-p$ free coefficients. Consequently, the optimization of the stability polynomials of the second- and third-order methods can in this case be performed using the monomial representation

$$P_{p;E}(z;\boldsymbol{\alpha}) = \sum_{j=0}^{p} \frac{1}{j!} z^j + \sum_{j=p+1}^{E} \alpha_j z^j \tag{4.9}$$

which is agnostic of the concrete form of the Butcher array.

If one would employ the unrestricted optimization of the stability polynomial in monomial form (4.5) one obtains $E-4$ monomial coefficients $\alpha_j$ which need to be mapped in some way onto the $E-5$ Butcher array coefficients $a_i$. The optimal way of doing this is unclear, and it is observed that setting the Butcher array coefficients $a_i$ to a least-squares solution leads to significantly inferior results compared to an optimization of the $\gamma_k$ using the tableau-related form (4.7).

*4.3. Loss of Optimality*

Due to the design choices of the new fourth-order P-ERK schemes discussed in Section 3.2 we face an inevitable loss of optimality in the stability polynomial optimization, as the in principle free monomial coefficient $\alpha_5$ is determined by the other coefficients $\alpha_j, j = 6, \dots, E$. Consequently, the stability polynomials which we use to construct the P-ERK methods will yield a smaller timestep than the ones with completely free coefficients $\alpha_j, j = 4, \dots, E$. In this section, we seek to quantify this loss of optimality by comparing the maximum admissible timestep $\Delta t_{\text{PERK}}$ of the P-ERK stability polynomial to the maximum admissible timestep $\Delta t$ of a 'free' fourth-order accurate stability polynomial for a concrete example.

To this end, we optimize stability polynomials of both kinds with degrees $E = \{5, 6, \dots, 16\}$ for the spectrum of the isentropic vortex advection testcase [39, 55]. The compressible Euler equations are discretized using the discontinuous Galerkin spectral element method (DGSEM) [56, 57] with Harten-Lax-Van Leer-Contact (HLLC) numerical flux [58] and local polynomials of degree $k = 3$. Throughout this work, we use `Trixi.jl` [59–61], a high-order DG code written in `Julia` for the numerical simulation of conservation laws. `Trixi.jl` supports algorithmic differentiation (AD), which allows exact and efficient computation of the Jacobian $J(\boldsymbol{U}_0) = \frac{\partial \boldsymbol{F}}{\partial \boldsymbol{U}}|_{\boldsymbol{U}_0}$, cf. (2.6), from which the spectrum $\boldsymbol{\sigma}(J) = \{\lambda_1, \dots, \lambda_M\}$ is computed. The spectrum of this setup on a $8 \times 8$ grid is shown in Fig. 1 together with the domains of absolute stability

$$\mathcal{S}_{4;E} \coloneqq \{z \in \mathbb{C} \,:\, |P_{4;E}(z)| \leq 1\} \tag{4.10}$$

of an unconstrained fourth-order and P-ERK stability polynomial. The optimal timestep $\Delta t^{(6)}$ of the sixth-degree, fourth-order unconstrained stability polynomial is used to scale the spectrum $\boldsymbol{\lambda}$, i.e., we plot $z_i = \Delta t^{(6)} \cdot \lambda_i$, $i = 1, \dots, M$. In addition, we plot the domains of absolute stability $\mathcal{S}$ for the unconstrained and P-ERK stability polynomials. Clearly, the P-ERK stability polynomial yields a smaller domain of absolute stability than the unconstrained one resulting in a smaller maximum admissible timestep $\Delta t^{(6)}_{\text{P-ERK}}$.

To assess the loss of optimality quantitatively, we present the absolute values $\Delta t$ and $\Delta t_{\text{PERK}}$ of the maximum admissible timesteps in Fig. 2a and the ratio ${}^{\Delta t_{\text{P-ERK}}}/_{\Delta t}$ in Fig. 2b. As expected, for all polynomial degrees $E$ the maximum admissible timesteps $\Delta t_{\text{P-ERK}}$ of the P-ERK stability polynomials are smaller than the maximum admissible timesteps $\Delta t$ of the free stability polynomials. However, already for the $E = 6$ P-ERK polynomial with one free coefficient, approximately 85% of the maximum admissible timestep $\Delta t$ of the two-coefficient P-ERK stability polynomial can be achieved. This ratio increases with $E$ up to 96.5% for the $E = 16$ P-ERK stability polynomial. Clearly, the price we pay for the design choices of the new fourth-order P-ERK methods is comparatively small and is a justifiable loss in optimality of the stability polynomials.

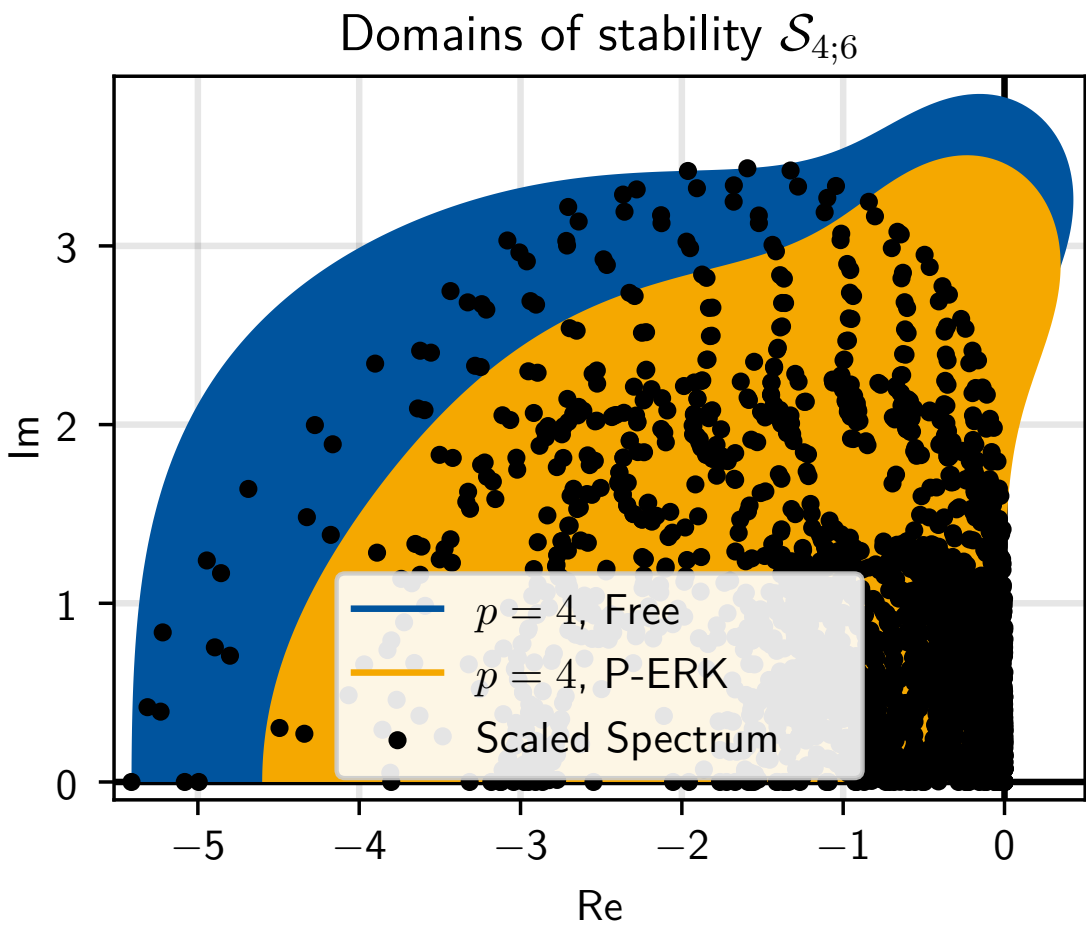

Figure 1: Scaled spectrum corresponding to a DGSEM discretization of the isentropic vortex advection (see description of the setup in Section 4.3). The spectrum $\boldsymbol{\lambda} \in \boldsymbol{\sigma}\left(\frac{\partial \boldsymbol{F}}{\partial \boldsymbol{U}}\right)$ is scaled by the optimal timestep $\Delta t^{(6)}$ of the sixth-degree, fourth order unconstrained stability polynomial, i.e., $z_m = \Delta t^{(6)} \lambda_m$, $m = 1, \dots, M$. The domains of absolute stability of the unconstrained and P-ERK stability polynomials are shown in blue and orange, respectively. Note that due to symmetry around the real axis we show only the second quadrant of the complex plane.

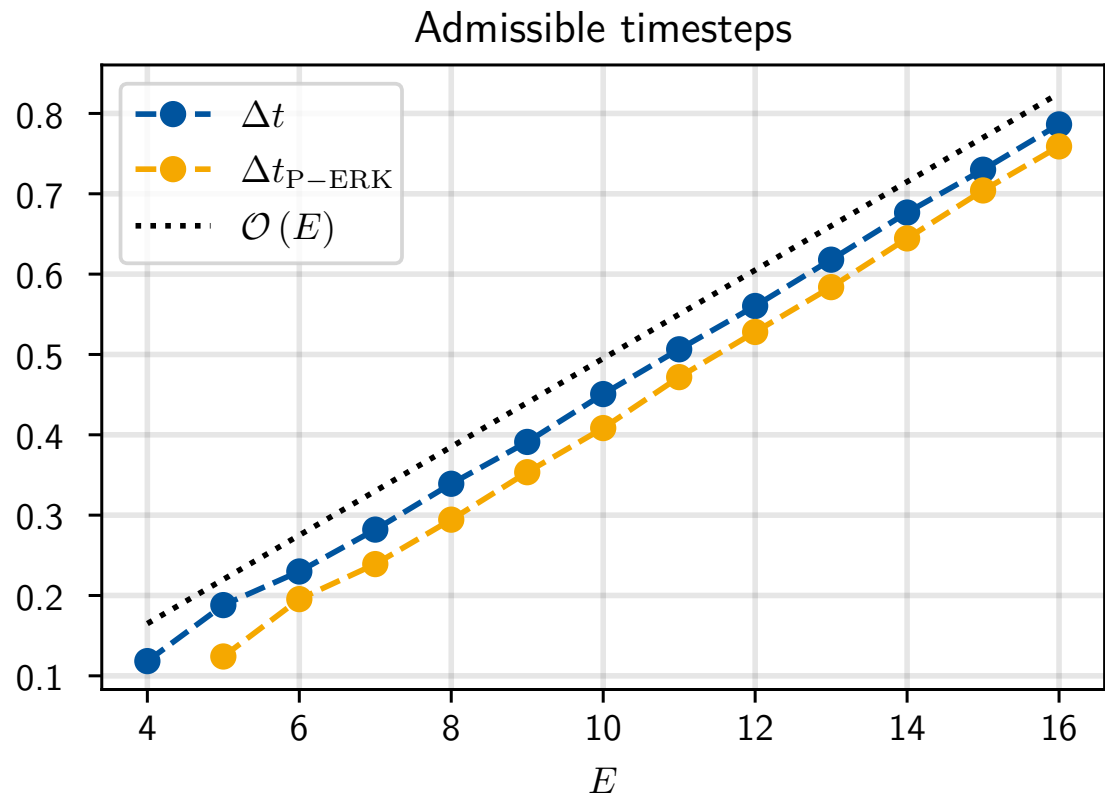

(a) Absolute values of the maximum admissible timesteps $\Delta t$ and $\Delta t_{\text{PERK}}$. Note the linear increase of the stable timesteps with number of stage evaluations $E$ (dotted black line).

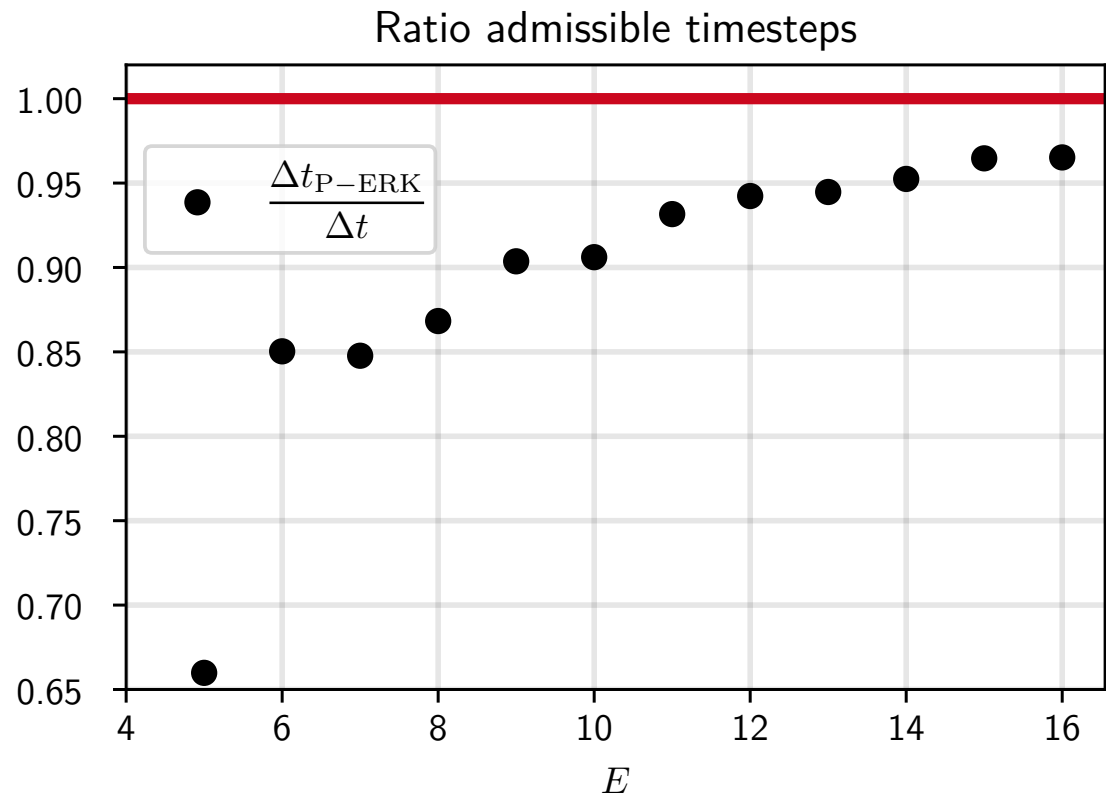

(b) Realizable optimality of the P-ERK stability polynomials compared to an unconstrained fourth-order stability polynomial.

Figure 2: Comparison of absolute Fig. 2a and relative Fig. 2b maximum admissible timesteps $\Delta t_{\text{PERK}}, \Delta t$ for the spectrum of a DGSEM discretization of the isentropic vortex advection testcase.

### *4.4. Internal Stability*

It is well-known that internal stability, i.e., the propagation of round-off errors is of considerable concern for many-stage methods [52, 62–65]. This issue becomes also relevant for the P-ERK methods as the number of stage evaluations $E$ increases, as displayed in Fig. 3. The exponential increase in $\widetilde{\mathcal{M}}$ with the number of stage evaluations $E$ showcases the vulnerability of the methods with many stage evaluations to round-off errors. In the presence of round-off errors due to limited floating point precision, the approximations $\boldsymbol{U}_n$ and stages $\boldsymbol{K}_i$ cf. (2.2) are perturbed by round-off errors $\boldsymbol{e}_i$ [52]. These perturbed quantities $\widetilde{\boldsymbol{U}}_n, \widetilde{\boldsymbol{K}}_i$ are computed in practice as [52]:

$$\widetilde{\boldsymbol{K}}_i = \boldsymbol{F}^{(r)}\left(t_n + c_i\Delta t, \widetilde{\boldsymbol{U}}_n + \Delta t\sum_{j=1}^{S} a_{i,j}\widetilde{\boldsymbol{K}}_j + \boldsymbol{e}_i\right), \quad i = 1, \dots, S \tag{4.11a}$$

$$\widetilde{\boldsymbol{U}}_{n+1} = \widetilde{\boldsymbol{U}}_n + \Delta t\sum_{i=1}^{S} b_i\widetilde{\boldsymbol{K}}_i + \boldsymbol{e}_{S+1}\,. \tag{4.11b}$$

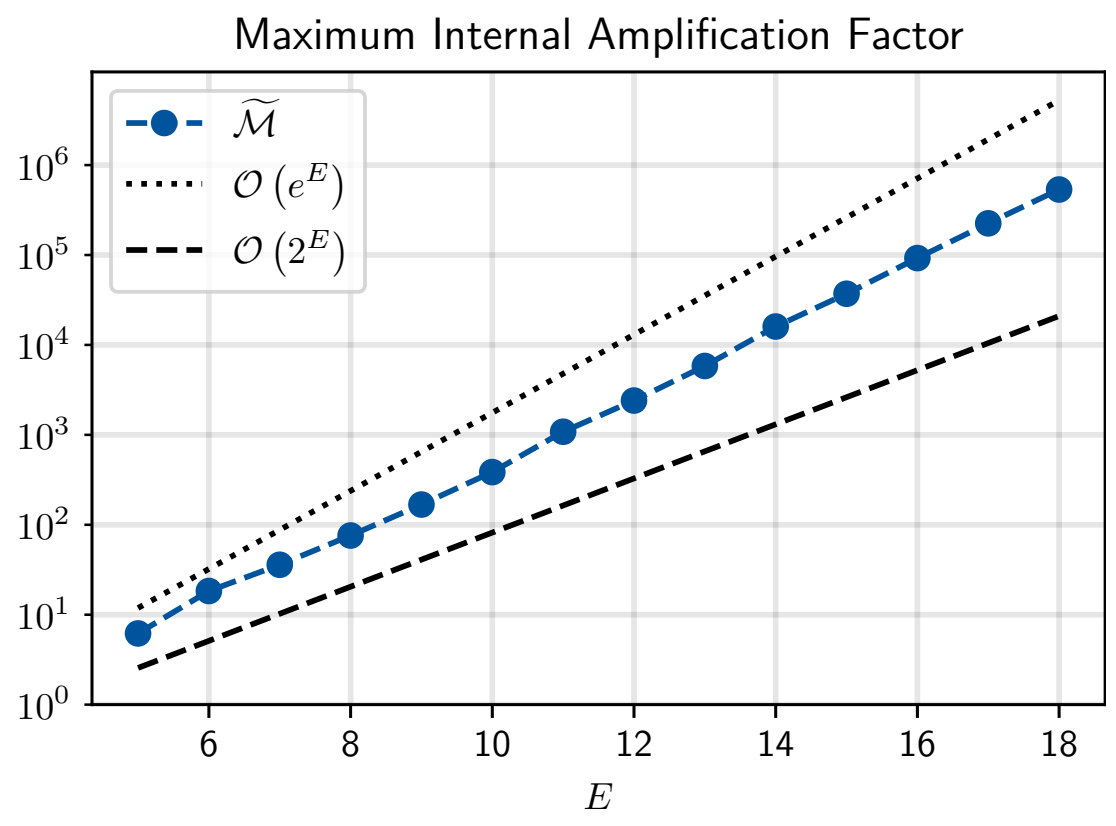

Figure 3: Maximum internal amplification factor $\widetilde{\mathcal{M}}$ (4.14) for $E = \{5, \dots, 18\}$ stage-evaluation fourth-order P-ERK methods which have been optimized for the spectrum of a DGSEM discretization of the isentropic vortex advection testcase. The depicted values $\widetilde{\mathcal{M}}$ are obtained for the uniform timesteps $c_i = 1.0, i = 2, \dots, E-3$.

In [52] it was shown that the defect from computed, perturbed solution $\widetilde{\boldsymbol{U}}_n$ to the exact solution $\boldsymbol{U}(t_n)$ can be estimated for linear ODEs as

$$\left\|\widetilde{\boldsymbol{U}}_{n+1} - \boldsymbol{U}(t_{n+1})\right\| \leq \left\|\widetilde{\boldsymbol{U}}_n - \boldsymbol{U}(t_n)\right\| + \sum_{i=1}^{S+1}\|Q_i(z)\|\,\|\boldsymbol{e}_i\| + \mathcal{O}\left(\Delta t^{p+1}\right). \tag{4.12}$$

In addition to the usual truncation error $\mathcal{O}\left(\Delta t^{p+1}\right)$ there is a term corresponding to the amplification of round-off errors $\sum_k \|Q_k(z)\|\|\boldsymbol{e}_k\|$. Internal stability becomes a concern when both are of similar magnitude [52].

In the estimate of the error of the perturbed iterates (4.12) $Q_j(z)$ denote the *internal stability polynomials* [63, 52] which for a method in Butcher form are computed as [52]

$$\boldsymbol{Q}(z; \boldsymbol{b}, A) \coloneqq z\boldsymbol{b}^T\left(I - zA\right)^{-1} \quad \in \mathbb{C}^{1\times S}, \tag{4.13}$$

where $I \in \mathbb{R}^{S\times S}$ denotes the identity matrix. For explicit methods there is no initial round-off error $\boldsymbol{e}_1$ [52] and it suffices to consider the internal stability polynomials starting from the second stage $(i = 2)$. To estimate the potential of round-off error amplification it is customary to investigate

$$\widetilde{\mathcal{M}}(\boldsymbol{b}, A) \coloneqq \max_{z\in\{\Delta t\lambda_m\}}\sum_{i=2}^{S}|Q_i(z; \boldsymbol{b}, A)| \tag{4.14}$$

which is a variant of the *maximum internal amplification factor* proposed in [52]. Although (4.12) is only valid for linear ODEs it is a useful heuristic for the potential of round-off error amplification also for nonlinear ODEs [52, 65].

For a given spectrum $\boldsymbol{\sigma} = \{\lambda_m\}_{m=1,\dots,M}$ and correspondingly optimized stability polynomial $P_{4;E}(z; \boldsymbol{\gamma})$ one can compute $\widetilde{\mathcal{M}}$ for different weight vectors $\boldsymbol{b}$ and Butcher arrays $A$. We seek to keep the simple form of the weight vector $\boldsymbol{b} = \begin{pmatrix}0 & \dots & 0 & 0.5 & 0.5\end{pmatrix}$ but test different Butcher arrays $A$. In particular, we examine $\widetilde{\mathcal{M}}$ for different choices of the free abscissae $c_i, i = 2, \dots, S-3$ which propagate into $A$, cf. (3.7). In principle, optimization of $\widetilde{\mathcal{M}}$ over the $c_i$ is possible, but due to the highly nonlinear objective (4.14) involving inversion of a matrix (4.13) this is a challenging task with little hope for success [65].

| Abscissae | $\widetilde{\mathcal{M}}$ | CFL |
|---|---|---|
| Linear Increasing $c_i = \frac{i}{E-4}$ | $1.35809 \cdot 10^6$ | 0.87 |
| Constant $c_i = 1.0$ | $5.34877 \cdot 10^5$ | 0.90 |

Table 1: Maximum internal amplification factor $\widetilde{\mathcal{M}}$ (4.14) and CFL number for $E = 18$ P-ERK methods with different choices of the abscissae $c_i$.

Instead, we tried linear increasing $c_i = {}^{i}/_{E-4}$ as in [1, 2] and constant $c_i = c_{\text{Const}}$ abscissae. Out of these two, constant timesteps with $c_{\text{Const}} = 1.0$ yielded the smaller amplification factor $\widetilde{\mathcal{M}}$ and are thus used in the following. This can be intuitively explained since this puts maximum weight on the first stage $\boldsymbol{K}^{(1)}$ which in the context of one timestep $t_n \to t_{n+1}$ is the least-perturbed stage. Furthermore, allowing for maximum abscissae $c_i \in (0, 1]$ reduces the risk of negative entries $a_{i,1}$ (cf. (3.7)), which correspond to downwinding the corresponding stage [40]. To illustrate the influence of the choice of the abscissae $c_i$ on the internal stability, in Table 1 we present the internal amplification factors $\widetilde{\mathcal{M}}$ and CFL number for $E = 18$ P-ERK methods. The stability polynomials have been optimized for the spectrum of the isentropic vortex advection testcase shown in Fig. 1 that is also employed in convergence study Section 5.4.1.

## 5. Validation

Before turning to applications we compare the combined P-ERK schemes to their composing standalone schemes with regard to errors and conservation of linear invariants. Additionally, we demonstrate linear stability and fourth-order convergence.

### *5.1. Linear/Absolute Stability*

To showcase that linear stability of each individual optimized Runge-Kutta method for its corresponding spectrum suffices for overall linear stability, we consider the 1D advection equation with smooth initial data

$$u_t + u_x = 0, \tag{5.1a}$$

$$u_0(t_0 = 0, x) = 1 + \frac{1}{2}\sin(\pi x) \tag{5.1b}$$

on $\Omega = (-1, 1)$ equipped with periodic boundaries. The advection equation is discretized using the DGSEM [56, 57] with local polynomials of degree $k = 3$ using `Trixi.jl` [59–61], a high-order DGSEM code written in native Julia. The mesh is discretized using 192 cells, which are refined by a factor two in the center of the domain $x \in [-0.5, 0.5]$. This results in a non-uniform grid with cell sizes $\Delta x^{(1)} = 2^{-7}$ and $\Delta x^{(2)} = 2^{-6}$.

Since the equation is linear, the semidiscretization (2.1a) may be written as

$$\boldsymbol{U}'(t) = \begin{pmatrix}\boldsymbol{U}^{(1)}(t)\\ \boldsymbol{U}^{(2)}(t)\end{pmatrix}' = \begin{pmatrix}A^{(1)} & 0\\ 0 & A^{(2)}\end{pmatrix}\boldsymbol{U}'(t)\,, \tag{5.2}$$

where the partitioning is performed according to the local cell size. We apply a two-level ($R = 2$) P-ERK scheme to the system (5.2). In this case, we optimize $E = \{10, 16\}$ schemes to the spectrum of a uniformly discretized advection equation. The admissible timestep $\Delta t$ of the $E = 16$ stage scheme is roughly twice as large as the timestep of the $E = 10$ stage scheme. We thus use the $E = 16$ stage scheme in the refined part of the domain, while the $E = 10$ stage scheme is used to integrate the coarse part. Applying the two-level P-ERK scheme to (5.1) yields the fully discrete system

$$\boldsymbol{U}_{n+1} = \underbrace{\begin{pmatrix} D^{(1)}\left(\Delta t, A^{(1)}\right) & 0 \\ 0 & D^{(2)}\left(\Delta t, A^{(2)}\right) \end{pmatrix}}_{=:D} \boldsymbol{U}_n \, , \tag{5.3}$$

where $D$ is equivalent to the matrix-valued stability function $P(Z)$ (2.12) with $R = 2$ and specified mask matrices $I^{(1)}, I^{(2)}$. As (5.3) is a constant-coefficient linear time-invariant system, the stability of the fully discrete system is entirely determined by the eigenvalues of the matrix $D$ which are plotted in Fig. 4. As the spectral radius is $\rho(D) \approx 1$ (up to machine precision), we conclude linear stability [38, Chapter 8.6], which is confirmed by long-running simulations up to $t_f = 100$ corresponding to 7672 timesteps.

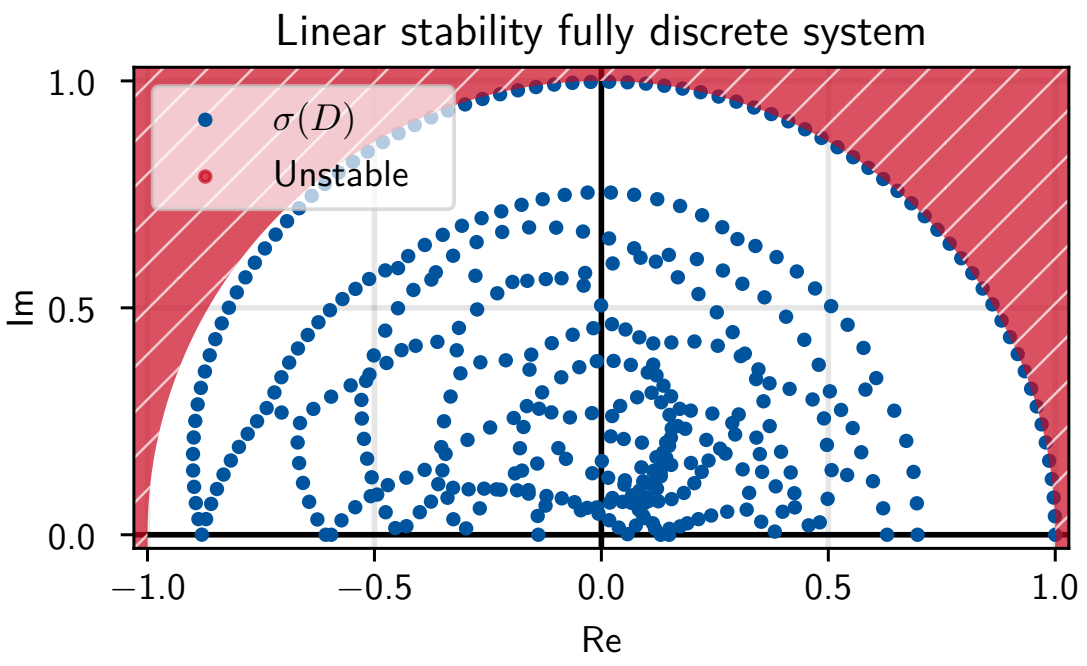

Figure 4: Spectrum $\boldsymbol{\sigma}(D)$ of the fully discrete matrix $D$ (5.3) of the P-ERK$_{4;\{10,16\}}$ scheme applied to (5.2). The unstable region with $|\lambda| > 1$ is shaded in red.

### *5.2. Error Comparison*

In addition to convergence studies presented in Section 5.4 that demonstrate the asymptotic convergence rate of a combined P-ERK scheme, we compare the errors of a family of P-ERK schemes to the errors of the individual schemes, i.e., the error constants. This information is essential for assessing the effectiveness of the P-ERK schemes since one is ultimately interested in accuracy per computational cost. This error comparison corresponds to a check for internal consistency of the involved schemes, as for internally inconsistent schemes significant errors at the interfaces between different schemes are expected [18, 21].

For the remaining tests in this section we consider again the classic isentropic vortex testcase [39, 55] with the same parameters as used in [1, 2] on domain $\Omega = [0, 10]^2$. For simulation of the compressible Euler equations we again use `Trixi.jl` [59–61]. The fluxes are approximated with the HLLC numerical flux [58]. The domain $\Omega$ is initially discretized with 64 elements per direction on which the solution is represented using polynomials with degree $k = 3$ for an overall fourth-order accurate method. The mesh is dynamically refined around the vortex (cf. Fig. 5a) based on the distance to the center of the vortex with two additional levels of refinement, i.e., the minimal grid size is in this case $h = \frac{10}{256}$. The mesh at final time $t_f = 20$ is displayed in Fig. 5b. We construct stability polynomials of degrees $E = \{6, 10, 16\}$ where the maximum admissible timestep roughly doubles for each increase in degree. The density errors at final time $t_f = 20$

$$e_\rho(x, y) \coloneqq \rho(t_f, x, y) - \rho_h(t_f, x, y) \tag{5.4}$$

are tabulated in Table 2 for each individual optimized Runge-Kutta method and the composed P-ERK scheme. Both $L^\infty$ and domain-normalized $L^1$ errors for the joint P-ERK schemes are of

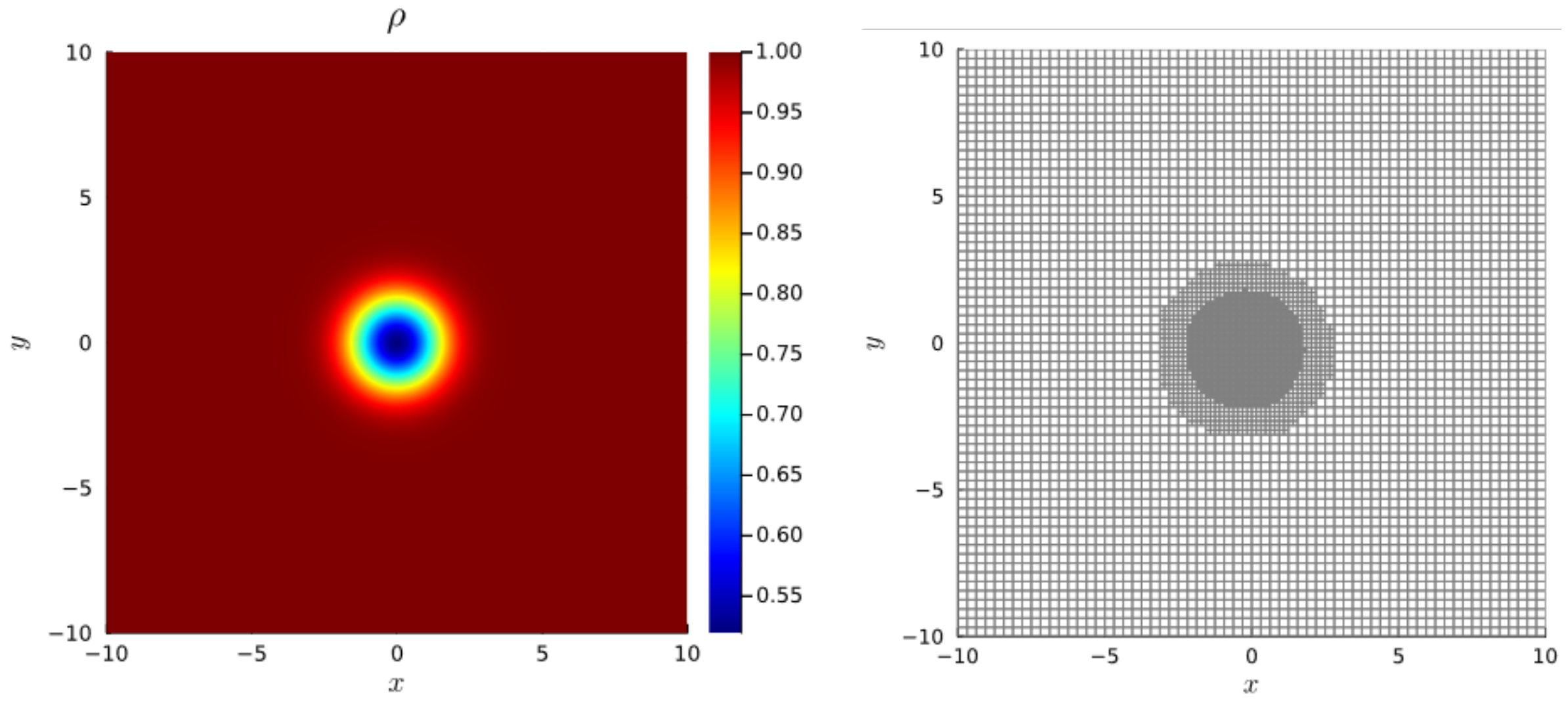

(a) Density $\rho$ at final time $t_f = 20$.

(b) Adaptively refined three-level mesh at final time $t_f = 20$.

Figure 5: Density $\rho$ and adaptively refined mesh at final time $t_f = 20$ corresponding to one pass through the domain $\Omega = [0, 10]^2$ for the isentropic vortex advection testcase [39, 55].

| Method | $\frac{1}{\lvert\Omega\rvert}\lVert e_\rho(x,y)\rVert_{L^1(\Omega)}$ | $\lVert e_\rho(x,y)\rVert_{L^\infty(\Omega)}$ |
|---|---|---|
| P-ERK$_{4;6}$ | $3.63200009 \cdot 10^{-7}$ | $1.14969098 \cdot 10^{-5}$ |
| P-ERK$_{4;10}$ | $3.63200072 \cdot 10^{-7}$ | $1.14968812 \cdot 10^{-5}$ |
| P-ERK$_{4;16}$ | $3.63199965 \cdot 10^{-7}$ | $1.14965181 \cdot 10^{-5}$ |
| P-ERK$_{4;\{6,10,16\}}$ | $3.63155882 \cdot 10^{-7}$ | $1.14988553 \cdot 10^{-5}$ |

Table 2: Comparison of density-errors for the isentropic vortex advection testcase after one pass through the domain $\Omega = [0, 10]^2$ at $t_f = 20$. The compressible Euler equations are discretized using the DGSEM with HLLC flux [58] on a dynamically refined mesh with cell sizes ranging from $\Delta x^{(1)} = {}^{10}/_{64}$ to $\Delta x^{(3)} = {}^{2}/_{256}$. The local polynomials are of degree $k = 3$ to enable an overall fourth-order accurate method.

similar size as the errors of the individual schemes. Especially the fact that there is no increase in $L^\infty$ norm shows that the P-ERK schemes can be used in practical simulations, despite their troubling nonlinear stability properties [43].

This is noteworthy as for instance in [18] it is shown that at the interfaces of different schemes order reduction [66] may be encountered. As argued in [18], the interfaces between methods act like additional time-dependent internal boundary conditions, thus potentially suffering from effects also encountered for non-periodic problems [66]. As shown in Section 5.4 fourth-order pointwise convergence is achieved despite the discontinuities of methods at the interfaces.

### *5.3. Conservation of Linear Invariants*

We seek to demonstrate that the P-ERK schemes conserve linear invariants such as mass, momentum, and energy according to Section 2.4. To this end, we compare the newly derived fourth-order P-ERK schemes to a well-established fourth-order Runge-Kutta method. In particular, we consider the 10-stage, fourth-order low-storage method by Ketcheson [67] which is implemented in the `OrdinaryDiffEq.jl` package [68]. As for the comparison of the density-errors of the P-ERK schemes to standalone schemes, the local polynomials are of degree $k = 3$ to yield an overall fourth-order scheme. The conservation error

$$e_u^{\text{Cons}}(t_f) \coloneqq \left| \int_\Omega u(t_0, \boldsymbol{x}) \mathrm{d}\boldsymbol{x} - \int_\Omega u(t_f, \boldsymbol{x}) \mathrm{d}\boldsymbol{x} \right| \tag{5.5}$$

| | P-ERK$_{4;\{6,10,16\}}$ | SSP$_{4;10}$ [67] |
|---|---|---|
| $e_{\rho}^{\text{Cons}}(t_f)$ | $3.78 \cdot 10^{-13}$ | $2.11 \cdot 10^{-13}$ |
| $e_{\rho v_x}^{\text{Cons}}(t_f)$ | $5.70 \cdot 10^{-14}$ | $1.26 \cdot 10^{-13}$ |
| $e_{\rho v_y}^{\text{Cons}}(t_f)$ | $6.53 \cdot 10^{-14}$ | $1.29 \cdot 10^{-13}$ |
| $e_{\rho e}^{\text{Cons}}(t_f)$ | $9.40 \cdot 10^{-13}$ | $3.66 \cdot 10^{-12}$ |

Table 3: Conservation-errors for isentropic vortex advection testcase after one pass through the domain at $t_f = 20$ for the fourth-order P-ERK scheme and the ten-stage, fourth-order SSP method [67].

| CFL | $\frac{1}{\lvert\Omega\rvert}\lVert e_\rho(x,y)\rVert_{L^1(\Omega)}$ | EOC | $\lVert e_\rho(x,y)\rVert_{L^\infty(\Omega)}$ | EOC |
|---|---|---|---|---|
| 1 | $3.37848279 \cdot 10^{-10}$ | | $1.30746968 \cdot 10^{-6}$ | |
| $1/2$ | $2.12288265 \cdot 10^{-11}$ | 3.99 | $4.23331622 \cdot 10^{-8}$ | 4.95 |
| $1/4$ | $6.47315254 \cdot 10^{-12}$ | 1.71 | $1.96120620 \cdot 10^{-9}$ | 4.43 |

Table 4: Convergence study for the isentropic vortex advection testcase by means of the density-errors after one pass through the domain at $t_f = 20$ for the fourth-order P-ERK scheme.

is computed for mass, momentum and energy after one pass of the vortex through the domain and tabulated in Table 3. Clearly, the errors of the $p = 4, E = \{6, 10, 16\}$ P-ERK family are on par with the errors of the SSP scheme. We remark that these conservation violations are caused by floating point induced errors like round-off, as the spatial discretization is conservative by construction.

### *5.4. Convergence Study*

We begin this section by performing a convergence study inspired by [1, 2] by means of the isentropic vortex advection testcase [39, 55]. Due to the fourth-order accuracy of the P-ERK schemes presented here and the relatively strict CFL condition for high degree local polynomials, only very few reductions in timestep $\Delta t$ may be performed until the spatial errors dominate.

To showcase fourth-order convergence in the asymptotic regime we thus resort to an ODE system with no underlying spatial discretization. To this end we consider the well-known Lotka-Volterra predator-prey system [69, 70] for which recently an analytical solution has been derived [70]. The Lotka-Volterra predator-prey system is a nonlinear, naturally partitioned ODE system and thus well-suited for a convergence study which we present in Section 5.4.2.

#### *5.4.1. Isentropic Vortex Advection*

As in [1, 2] we employ local polynomials of degree $k = 6$ for the convergence tests for the fourth-order P-ERK scheme to minimize spatial errors. As for the previous tests we use a base discretization composed of $64^2$ elements and refine the mesh around the vortex with two additional levels of refinement. We construct a family of $E = \{5, 9, 15\}$ fourth-order schemes such that the admissible timestep doubles for each family member.

Despite this high spatial accuracy, the errors of the P-ERK schemes become quickly dominated by spatial errors again when decreasing the timestep. This is demonstrated in Table 4. We clearly have the expected fourth-order convergence in $L^1$-norm when going from CFL = 1 to CFL = $1/2$. For the $L^\infty$-norm we observe a slightly higher order of convergence, which we attribute to the fact that we are not yet in the asymptotic regime. At CFL = $1/4$ the $L^1$ error is already dominated by spatial errors, which is reflected in the reduced order of convergence. In particular, the smallest $L^1$ error of the $p = 3$ P-ERK scheme observed in the convergence study performed in [43] is equal to $6.26 \cdot 10^{-12}$ which is almost identical to the error of the $p = 4$ P-ERK scheme.

*5.4.2. Lotka-Volterra Predator-Prey System*

In order to not be limited by spatial discretization errors we consider a pure ODE problem. The Lotka-Volterra predator-prey system [69, 70] can be written in simplest form as

$$\dot{u}(t) = u(1 - v) \tag{5.6a}$$

$$\dot{v}(t) = v(u - 1) \tag{5.6b}$$

where $u(t)$ denote the number of prey and $v(t)$ the number of predators. The system (5.6) is a simplification of the classic Lotka-Volterra system obtained by setting all model parameters to one, i.e., $\alpha = \beta = \gamma = \delta = 1$. This simplification allows for the derivation of an analytical, closed-form solution involving only a quadrature problem [71]. The Lotka-Volterra system is a naturally partitioned ODE system with $R = 2$ that is nonlinear in both right-hand-sides and thus well-suited to demonstrate fourth-order convergence of the P-ERK methods developed in this work. For solving the ODE we employ the $E = \{5, 9\}$ P-ERK schemes which were constructed for the isentropic vortex advection testcase considered above. We study the $L^\infty$ error

$$e_\infty \coloneqq \left\| \begin{pmatrix} u_h(t^\star) \\ v_h(t^\star) \end{pmatrix} - \begin{pmatrix} u(t^\star) \\ v(t^\star) \end{pmatrix} \right\|_\infty \tag{5.7}$$

and the average error

$$\bar{e} \coloneqq \left| \frac{u_h(t^\star) + v_h(t^\star)}{2} - \frac{u(t^\star) + v(t^\star)}{2} \right| \tag{5.8}$$

where the final time $t^\star$ marks the point where the two populations are equal [71], as shown in Fig. 6a. In (5.7) and (5.8) the numerical approximation is denoted with subscript $h$. Both errors are presented in Fig. 6b showcasing the expected fourth-order convergence. The $e_\infty$ error (5.7)

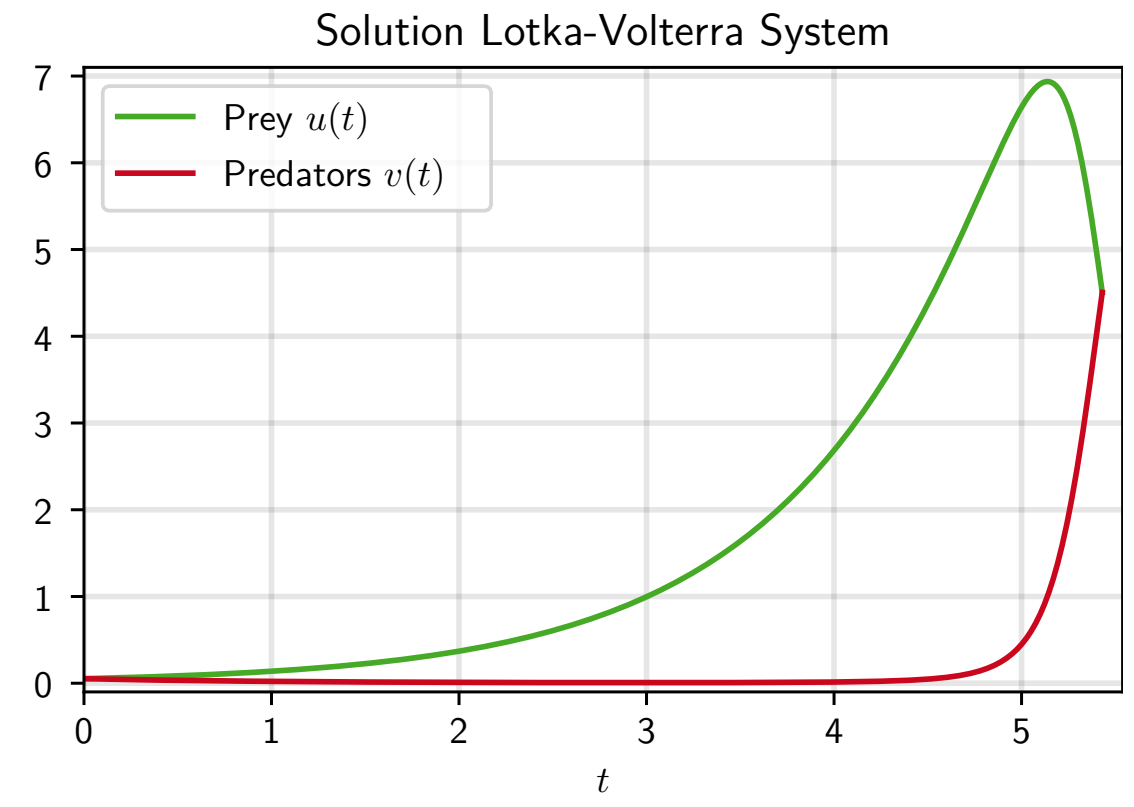

(a) Numerical solution to the Lotka-Volterra system with invariant $h = 2.0$ [71] over one period $(0, t^\star)$.

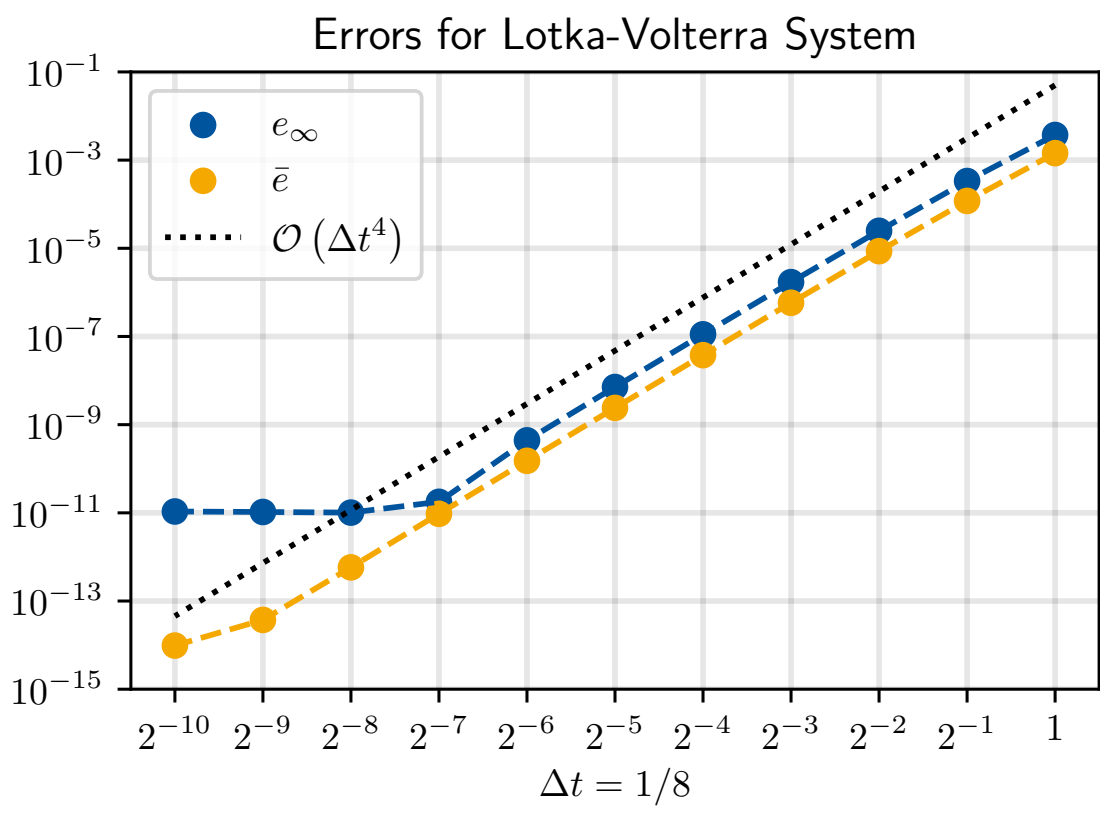

(b) Fourth-order convergence for the Lotka-Volterra predator-prey system (5.6). The largest timestep is given by $\Delta t = 1/8$ which is then halved until the smallest timestep $\Delta t = 2^{-13}$ is reached.

Figure 6: Solution and convergence study for the Lotka-Volterra predator-prey system solved with the $E = \{5, 9\}$ P-ERK scheme used also for the isentropic vortex testcase in Section 5.4.1.

starts to flatten at $\mathcal{O}(10^{-11})$ which is a consequence of the quadrature problem involved in the analytical solution [71]. In particular, the integrand becomes singular at the boundaries of the integral and is thus difficult to compute numerically. In contrast, for the averaged error $\bar{e}$ (5.8) the expected fourth-order convergence is observed almost up to machine precision $\mathcal{O}(10^{-15})$.

## 6. Applications

Before coming to applications we demonstrate the non-intrusive implementation of the partitioned time-integration for MoL codes in Section 6.1. Therein, we present one of the core

```julia
function apply_Jacobian!(du, mesh, equations, dg::DG, cache)
  @unpack inverse_Jacobian = cache.elements
  # Loop over all elements
  @threaded for element in eachelement(dg, cache)
    for j in eachnode(dg), i in eachnode(dg)
      factor = -inverse_Jacobian[i, j, element]
      for v in eachvariable(equations)
        du[v, i, j, element] *= factor
      end
    end
  end
  return nothing
end
```

Listing 1: 2D Reference-to-physical space transform for standard time integration schemes in `Trixi.jl` [59–61].

routines of the DG solver implemented in `Trixi.jl` for both the standard implementation and the required changes.

We consider four applications in this work representing different areas of applicability of the P-ERK schemes. First, we consider the linearized Euler equations in three spatial dimensions with non-uniform mean speed of sound $\bar{c} = \bar{c}(\boldsymbol{x})$ implying strongly varying characteristic speeds $\nu(\boldsymbol{x})$, cf. (3.2). The remaining three examples bear locally varying wave speeds due to non-uniform grids. Building upon our work in [43] we apply the fourth-order P-ERK schemes to adaptively refined meshes by means of the isentropic vortex advection testcase. Next, we consider the hyperbolic-parabolic flow around the SD7003 airfoil in the laminar regime, which is also employed as a testcase in [1, 2]. Finally, we apply the P-ERK schemes in a multi-physics setting by simulating the sound generation of two spinning, co-rotating vortices. In this case, two different P-ERK families are constructed, one for the nonlinear, compressible Euler equations and one for the linear acoustic perturbation equations.

### *6.1. Implementation in `Trixi.jl`*

Depending on the mesh type, there are up to four partitioning-indicating data structures to manage in `Trixi.jl`. For every mesh type, the `elements`, i.e., cells, need to be assigned to a certain partition $r$. For the tree-based meshes `TreeMesh`, `P4estMesh` [72] also `interfaces`, `boundaries`, and `mortars` need to be assigned to a certain partition $r$. The `boundaries` are simply assigned to the partition of the corresponding `element`. The `interfaces` and `mortars` between non-conforming, i.e., differently sized, elements [73] are assigned to the partition of the smaller element.

To highlight the non-intrusive nature of the P-ERK schemes, we present the back-transform from reference to physical element coordinates for the standard implementation in `Trixi.jl` in Listing 1 and the modification required for the P-ERK schemes in Listing 2. Clearly, no functional changes to the `apply_Jacobian!` function are required except for the restriction to a certain subset of elements $r$. This demonstrates the non-intrusive, easy-to-implement property of the P-ERK schemes for codes employing a MoL approach.

The timestep for every time integration scheme considered in this work CFL is based, i.e., a constant CFL number cf. (3.1) is set for the whole simulation. The construction of adaptive error-based timestep controllers as in [44] is left for future work.

### *6.2. 3D Linearized Euler Equations with non-uniform mean Speed of Sound*

Before coming to applications of the fourth-order P-ERK schemes on non-uniform meshes, we consider an example where locally varying characteristic speeds $\nu$ are induced through altering the spectral radii $\rho_i$ (cf. (3.2)), i.e., the signal speeds of the directional flux Jacobians. To this

```
function apply_Jacobian!(du, mesh, equations, dg::DG, cache,
                         elements_r::Vector{Int64}) # Additional argument
  @unpack inverse_Jacobian = cache.elements
  # Loop over elements of the r'th level
  @threaded for element in elements_r
    for j in eachnode(dg), i in eachnode(dg)
      factor = -inverse_Jacobian[i, j, element]
      for v in eachvariable(equations)
        du[v, i, j, element] *= factor
      end
    end
  end
  return nothing
end
```

Listing 2: 2D Reference-to-physical space transform for P-ERK schemes in `Trixi.jl` . In contrast to the standard transform shown in Listing 1, the P-ERK version takes an additional argument `elements_r` which contains the elements bundling the $k+1$ polynomial coefficients from $\boldsymbol{U}^{(r)}$, cf. (2.1a). Then, the loop over elements is restricted to the elements of the $r$'th level.

extend, we consider the linearized Euler equations in three dimensions with a non-uniform mean speed of sound $\bar{c}(\boldsymbol{x})$. The linearized Euler equations considered in this work are given by

$$\partial_t \underbrace{\begin{pmatrix} \rho' \\ v_1' \\ v_2' \\ v_3' \\ p' \end{pmatrix}}_{\boldsymbol{u}} + \partial_x \underbrace{\begin{pmatrix} \bar{\rho} v_1' + \bar{v}_1 \rho' \\ \bar{v}_1 v_1' + \frac{p'}{\bar{\rho}} \\ \bar{v}_1 v_2' \\ \bar{v}_1 v_3' \\ \bar{v}_1 p' + \bar{c}^2 \bar{\rho} v_1' \end{pmatrix}}_{=\boldsymbol{f}_1} + \partial_y \underbrace{\begin{pmatrix} \bar{\rho} v_2' + \bar{v}_2 \rho' \\ \bar{v}_2 v_1' \\ \bar{v}_2 v_2' + \frac{p'}{\bar{\rho}} \\ \bar{v}_2 v_3' \\ \bar{v}_2 p' + \bar{c}^2 \bar{\rho} v_2' \end{pmatrix}}_{=\boldsymbol{f}_2} + \partial_z \underbrace{\begin{pmatrix} \bar{\rho} v_3' + \bar{v}_3 \rho' \\ \bar{v}_3 v_1' \\ \bar{v}_3 v_2' \\ \bar{v}_3 v_3' + \frac{p'}{\bar{\rho}} \\ \bar{v}_3 p' + \bar{c}^2 \bar{\rho} v_3' \end{pmatrix}}_{=\boldsymbol{f}_3} = \boldsymbol{0} \tag{6.1}$$

where the overbar quantities $\overline{(\cdot)}$ denote the mean flow and the primed quantities $(\cdot)'$ the perturbations thereof. The spectral radii $\rho_i, i = 1, 2, 3$ of the directional flux Jacobians $J_i = \frac{\partial \boldsymbol{f}_i}{\partial \boldsymbol{u}}, i = 1, 2, 3$ are proportional to the speed of sound $\bar{c}$, which in this example is non-constant in space, i.e., $\bar{c} = \bar{c}(\boldsymbol{x})$. A physical interpretation of this scenario is a domain with spatially varying temperature $T(\boldsymbol{x})$, which relates to the speed of sound for an ideal gas via $\bar{c} = \sqrt{\gamma RT}$ with specific gas constant $R$ and heat capacity ratio $\gamma$. Recalling the estimate for the characteristic wave speeds for DG methods (3.2) this leads to a non-uniform timestep restriction (3.1) which motivates the use of the P-ERK schemes.

We consider the cube $\Omega = [0, 1]^3$ equipped with periodic boundaries and uniform discretization through $64^3$ elements. The spatial derivatives are discretized using the DGSEM using solution polynomials of degree $k = 3$ with Harten-Lax-Van Leer (HLL) two-wave numerical flux using the wave speed estimates by Davis [74, 75]. The mean flow is given by $\bar{\rho} = 1.0$, $\bar{v_i} = 1.0, i = 1, 2, 3$, and $\bar{c}$ is defined as $\bar{c}(\boldsymbol{x}) = 10 - 9d(\boldsymbol{x})$, where the distance function $d(\boldsymbol{x})$ is defined as

$$d(\boldsymbol{x}) = \begin{cases} \|\boldsymbol{x}\|_2 & \|\boldsymbol{x}\|_2 \leq 1 \\ 1 & \|\boldsymbol{x}\|_2 > 1 \end{cases}, \tag{6.2}$$

i.e., the speed of sound linearly decreases with radius $r = \|\boldsymbol{x}\|_2$ over the sphere with radius one from center to surface. This is illustrated in Fig. 7a where we show the subdomain $[0, 1]^3 \subset \Omega$. The non-uniform speed of sound is discretized analogously to the unsteady perturbations using the DGSEM, resulting in a system with more than 100 million degrees of freedom. The initial perturbations correspond to an acoustic wave with $v_i = \exp\left(-30\|\boldsymbol{x}\|_2^2\right), i = 1, 2, 3$, $\rho' = p' = -v_i$. A two-dimensional plot of the solution at $t = 0.3$ through the $z = 0$ plane is provided in Fig. 7b.

| Method | $\tau/\tau^\star$ | $N_{\text{RHS}}/N^\star_{\text{RHS}}$ |
|---|---|---|
| P-ERK$_{4;\{5,6,\dots,15\}}$ | 1.0 | 1.0 |
| P-ERK$_{4;15}$ | 1.41 | 1.75 |
| NDB$_{4;14}$ [77] | 1.60 | 1.80 |
| PKD$_{4;18}$ [78] | 2.56 | 2.82 |
| CK$_{4;5}$ [76] | 2.57 | 2.57 |
| RK$_{4;4}$ | 3.86 | 3.42 |

Table 5: Runtimes and number of scalar RHS evaluations (6.3) of different fourth-order integrators compared to the optimized $p = 4, E = \{5, 6, \dots, 15\}$ P-ERK integrator for the 3D linearized Euler equations with non-uniform mean speed of sound.

We construct ten P-ERK schemes with stage evaluations $E = \{5, 6, 7, 9, 10, 11, 12, 13, 14, 15\}$ corresponding to $\bar{c} = \{1, 2, \dots, 10\}$ which are applied according to the local speed of sound to keep the overall timestep constant. We compare the P-ERK schemes to different fourth-order methods implemented in `DifferentialEquations.jl` [68]. These are the five-stage, $2N$ storage method CK$_{4;5}$ designed for hyperbolic PDEs [76], the 14-stage low-storage method NDB$_{4;14}$ optimized for advection-dominated problems [77], and the 18-stage three-register method PKD$_{4;18}$ optimized for the spectral difference method for wave-propagation problems [78]. We also compare to the canonical four-stage fourth-order Runge-Kutta method RK$_{4;4}$. As this method has no degrees of freedom which can be optimized for the problem at hand, it is expected to perform poorly in comparison to the other methods. Nevertheless, we include this here, as in the previous publications on the P-ERK schemes [1, 2] also the non-optimized three-stage, third-order method by Shu and Osher was used as a reference case. In particular, the speedups reported in [1, 2] should be compared to the speedups we report here over the RK$_{4;4}$ method. Finally, we also run this example with the standalone optimized P-ERK$_{4;15}$ scheme.

The observed ratios of measured to optimal (i.e., P-ERK-Multi) runtime $\tau/\tau^\star$ and scalar (i.e., per degree of freedom (DoF)) RHS evaluations

$$N_{\text{RHS}} = N_t \cdot N_K \cdot M \tag{6.3}$$

are presented in Table 5. Here, $M$ is the dimension of the semidiscretization (1.2), $N_t$ the number of timesteps taken, and $N_K$ the number of stage evaluations per timestep. As expected, the multirate P-ERK schemes outperform the other methods in terms of runtime and number of scalar RHS evaluations, followed by the standalone optimized P-ERK$_{4;15}$ scheme. The NDB$_{4;14}$ scheme is performing comparably well, while the CK$_{4;5}$ and PKD$_{4;18}$ are roughly on-par in terms of performance. Unsurprisingly, the RK$_{4;4}$ scheme performs worst among the considered methods. Interestingly, RK$_{4;4}$ as implemented in `OrdinaryDiffEq.jl` [68] is disproportionally slow relative to the ratio of scalar RHS evaluations, which is not the case for every other scheme. This might be due to the fact that the RK$_{4;4}$ scheme is not implemented in an optimized, i.e., low-storage manner as the other methods.

By comparing the number of scalar RHS evaluations, one observes that the standalone P-ERK$_{4;15}$ scheme takes roughly 75% more evaluations than the multirate P-ERK schemes, but this manifests only in 41% higher runtime, i.e., only roughly 50% of the theoretical possible speedup is realized. This is due to fact that `Trixi.jl` is memory-bound for the standard DGSEM with a simple numerical flux and thus saved computations do not translate proportionally into runtime savings. In particular, data locality is in general lost for the partitioned approach, i.e., the coefficients in the $\boldsymbol{U}^{(r)}$ are not necessarily stored in a contiguous memory block. Also, the data structures containing `elements` and `interfaces` are accessed in a non-sequential manner which leads to cache misses and thus reduced performance. In contrast, for the standard methods everything is processed as present in memory which explains the non-ideal performance of the partitioned approach.

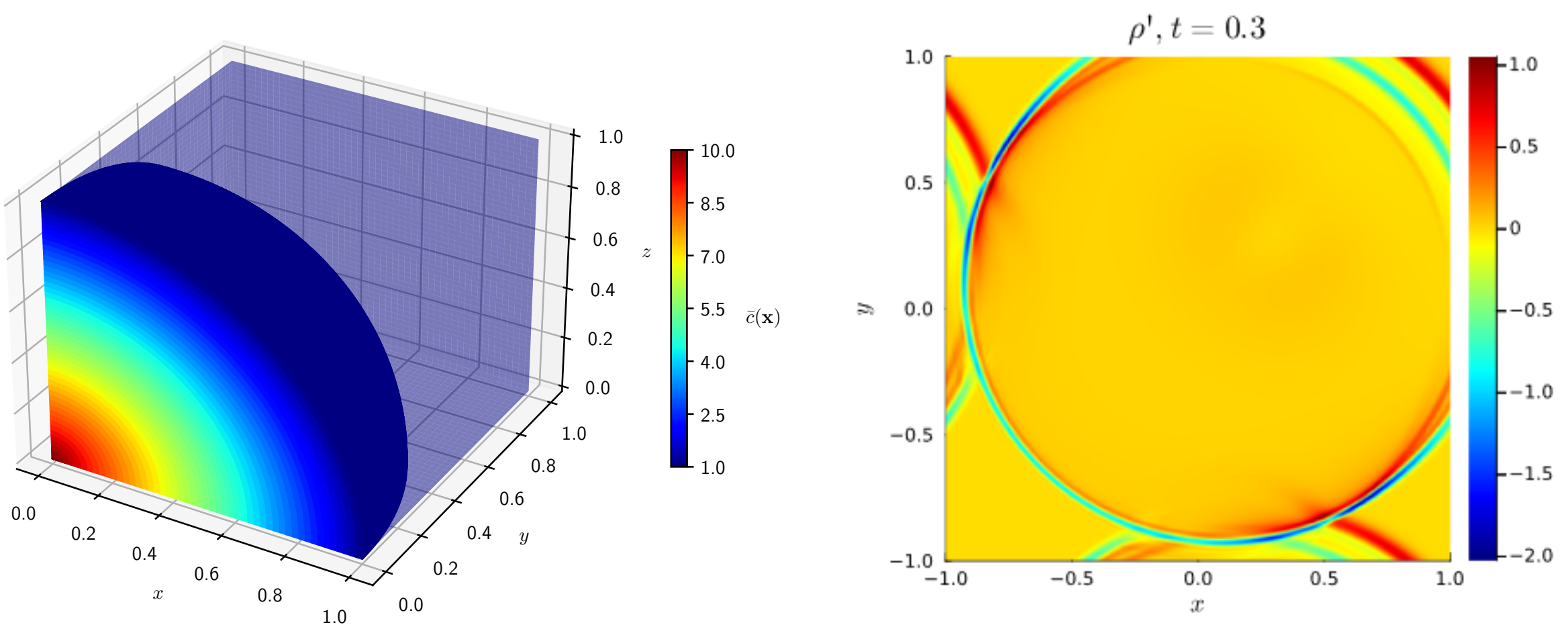

(a) Schematic of spatially varying speed of sound $\bar{c}(\boldsymbol{x})$ on the subdomain $[0,1]^3 \subset \Omega$.

(b) Planar view ($z = 0$) of the density perturbation $\rho'$ at $t = 0.3$.

Figure 7: Spatially varying speed of sound $\bar{c}(\boldsymbol{x})$ (Fig. 7a) and planar view of the density perturbation $\rho'$ at $t = 0.3$ for the 3D linearized Euler equations (Fig. 7b).

### *6.3. Isentropic Vortex Advection with Adaptive Mesh Refinement*

As a first application of the P-ERK schemes to non-uniform meshes we reconsider the isentropic vortex advection testcase employed already in Section 5. The parameters of the testcase are identical with the previous test, but we allow for an additional level of mesh refinement, resulting in a four-level grid. We have explored the application of P-ERK schemes to dynamically partitioned grids in [43] and found that the P-ERK schemes offer speedup compared to standard methods. Since for these problems often a significant share of cells lies on the finest level, the speedup is in general less pronounced compared to static non-uniform meshes used, e.g., for airfoil simulations.

For the isentropic vortex testcase we construct the four-member optimized $p = 4, E = \{5, 7, 11, 19\}$ P-ERK family and compare it to the standalone $p = 4, E = 19$ scheme. Additionally, we run the testcase for a collection of optimized fourth-order schemes implemented in `OrdinaryDiffEq.jl` [68]. These are the 14-stage low-storage method $\text{NDB}_{4;14}$ optimized for advection-dominated problems [77], the eight-stage, low-storage, low-dissipation, low-dispersion scheme $\text{TD}_{4;8}$ [79] optimized for DG methods applied to wave propagation problems, and the nine-stage method [80] optimized for spectral element discretizations of compressible fluid mechanics. Additionally, we also compare to the canonical four-stage fourth-order Runge-Kutta method $\text{RK}_{4;4}$.

The relative runtimes and scalar RHS evaluations are tabulated in Table 6. The multirate P-ERK scheme outperforms all other schemes, despite the additional overhead of the dynamic re-assignment of the partitioning data structures. The second fastest method is the standalone optimized $p = 4, E = 19$ scheme which is about 20% slower and requires 38% more scalar RHS evaluations. The optimized schemes available in `OrdinaryDiffEq.jl` perform roughly on-par and are at least 40% slower than the multirate P-ERK schemes. For this case, the canonical fourth-order Runge-Kutta method is surprisingly efficient, being only about 50% slower than the fastest standalone method. As for the previous example we observe that savings in scalar RHS evaluations do not translate exactly proportional into runtime savings.

### *6.4. Laminar Flow around SD7003 Airfoil*

As for the second- and third-order P-ERK schemes we apply the fourth-order P-ERK schemes to a laminar flow setup around the SD7003 airfoil. The flow parameters for this Navier-Stokes-

| Method | $\tau/\tau^\star$ | $N_{\text{RHS}}/N^\star_{\text{RHS}}$ |
|---|---|---|
| P-ERK$_{4;\{5,7,11,19\}}$ | 1.0 | 1.0 |
| P-ERK$_{4;19}$ | 1.20 | 1.38 |
| NDB$_{4;14}$ [77] | 1.50 | 1.82 |
| TD$_{4;8}$ [79] | 1.52 | 1.71 |
| RDPK$_{4;9}$ [80] | 1.40 | 1.63 |
| RK$_{4;4}$ | 1.79 | 2.08 |

Table 6: Runtimes and number of scalar RHS evaluations (6.3) of different fourth-order integrators compared to the optimized $p = 4, E = \{5, 7, 11, 19\}$ integrator for the isentropic vortex with four-level adaptive mesh.

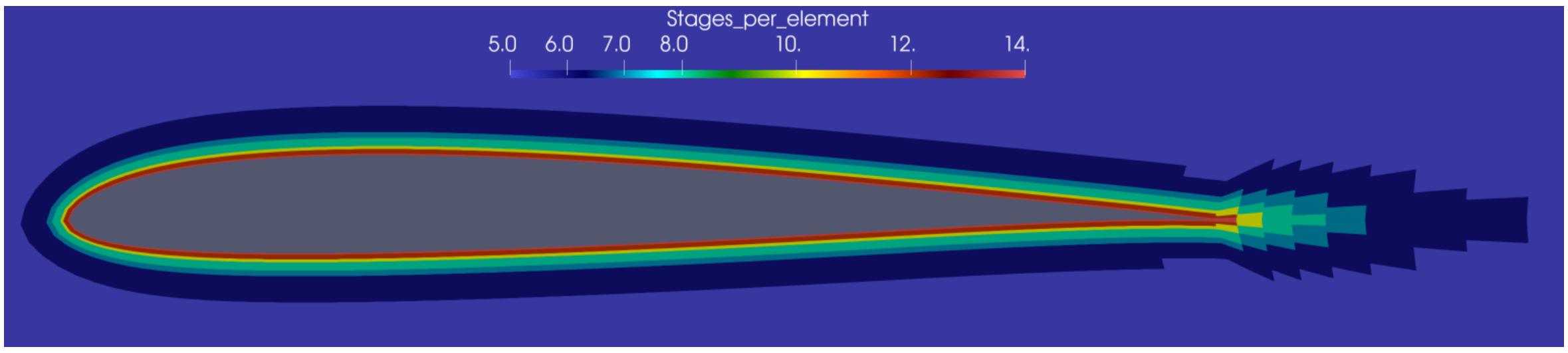

Figure 8: Method distribution for the SD7003 airfoil integration with the $E = \{5, 6, 7, 8, 10, 12, 14\}$ multirate P-ERK family.

Fourier simulation are $\text{Ma}_\infty = 0.2 = \frac{U_\infty}{a_\infty}$, $\text{Re}_\infty = 10^4 = \frac{\rho_\infty U_\infty c}{\mu_\infty}$ with airfoil chord length $c = 1.0$, angle of attack $\alpha = 4°$, non-dimensionalized speed of sound $a_\infty = 1.0$, pressure $p_\infty = 1.0$, and density $\rho_\infty = 1.4$ corresponding to isentropic exponent $\gamma = 1.4$. The viscosity is kept constant at $\mu = 2.8 \cdot 10^{-5}$ and the Prandtl number is set to $\text{Pr} = 0.72$. The reference state at the domain boundaries is weakly enforced (with no ramp-up) and the airfoil is modeled as a no-slip, adiabatic wall.

We employ a two-dimensional mesh consisting of 7605 straight-sided quadrilateral elements. This mesh (with second-order, i.e., curved-boundary elements) has been provided as electronic supplementary material to the second-order P-ERK paper [1]. The spatial derivatives are discretized using the DGSEM with $k = 3$ polynomials and HLLC numerical flux [58] with volume flux differencing [81] and entropy-stable volume flux [82]. The gradients required for the viscous fluxes are computed with the BR1 scheme [83, 84]. The $E = \{5, 6, 7, 8, 10, 12, 14\}$ P-ERK schemes are optimized for the free-stream state. The partitioning of the near-field grid is depicted in Fig. 8. We compare the multirate scheme against the P-ERK$_{4;12}$ scheme (which proved slightly more efficient than the 14-stage method), the NDB$_{4;14}$ and TD$_{4;8}$ schemes used also in Section 6.5, and the KCL$_{4;9}$ low-storage scheme [85] optimized for the compressible Navier-Stokes equations. The multirate P-ERK schemes require at most roughly half the number of scalar RHS evaluations compared to the other schemes. This leads to runtime savings between 35% for the standalone P-ERK scheme and 71% for TD$_{4;8}$. Interestingly, RK$_{4;4}$ is more efficient than KCL$_{4;9}$ which is optimized for the Navier-Stokes equations [85].

Prior to recording lift and drag coefficients, the simulation is run up to non-dimensional time $t_c := t\,\frac{U_\infty}{c} = 30$ to ensure all transients due to the unphysical initialization have decayed. A plot of the near-airfoil $y$-component of the velocity at $t = 30$ is shown in Fig. 9 which indicates the offstream of vortices from the sharp trailing edge which cause the lift and drag oscillations. Then, lift and drag coefficients are recorded over the $[30t_c, 35t_c]$ interval sampled with $\text{d}t_c = 5 \cdot 10^{-3}$ and reported in Table 7. The provided lift coefficient is solely due to the pressure difference:

$$C_L = C_{L,p} = \oint_{\partial\Omega} \frac{p\,\boldsymbol{n} \cdot \boldsymbol{t}^\perp}{0.5\rho_\infty U_\infty^2 L_\infty} \,\text{d}S \tag{6.4}$$

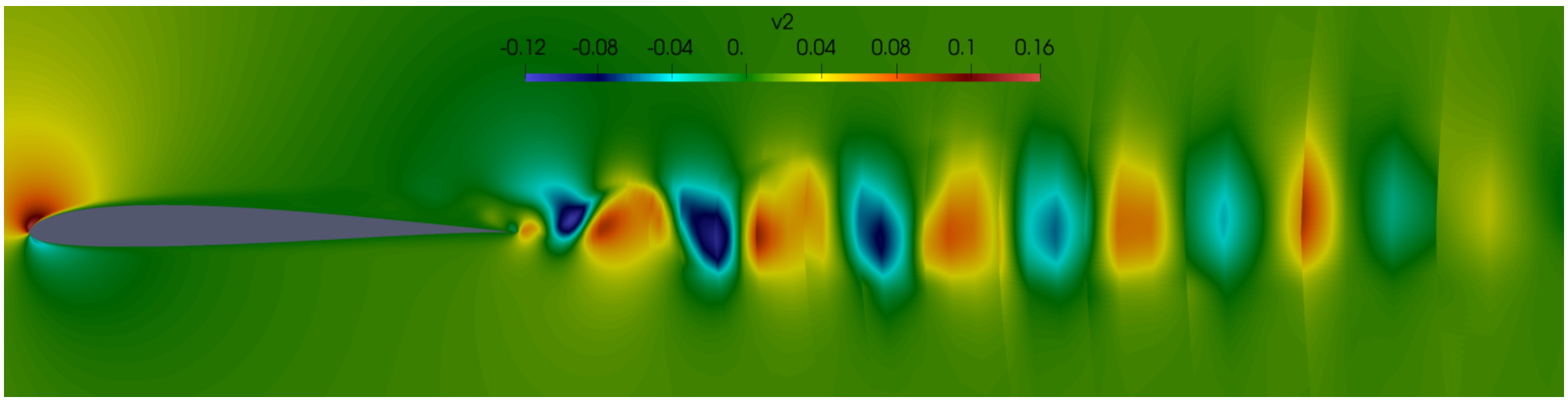

Figure 9: Vertical velocity $v_y$ at $t = 30$ for the SD7003 airfoil integration with the multirate P-ERK family.

| Method | $\tau/\tau^\star$ | $N_{\text{RHS}}/N^\star_{\text{RHS}}$ |
|---|---|---|
| P-ERK$_{4;\{5,6,\dots 14\}}$ | 1.0 | 1.0 |
| P-ERK$_{4;12}$ | 1.35 | 2.05 |
| NDB$_{4;14}$ [77] | 1.77 | 2.01 |
| TD$_{4;8}$ [79] | 1.71 | 1.90 |
| KCL$_{4;9}$ [85] | 2.82 | 3.10 |
| RK$_{4;4}$ | 2.28 | 2.33 |

Table 7: Runtimes and number of scalar RHS evaluations (6.3) of different fourth-order integrators compared to the optimized $p = 4, E = \{5, 6, 7, 8, 10, 12, 14\}$ integrator for the laminar flow around the SD70003 airfoil.

with $\boldsymbol{t}^\perp = \big(-\sin(\alpha), \cos(\alpha)\big)$ being the unit direction perpendicular to the free-stream flow and $\boldsymbol{n}$ as the fluid-element outward-pointing normal, i.e., the normal pointing into the surface. The drag coefficient incorporates also the viscous stress contribution:

$$C_D = C_{D,p} + C_{D,\mu} \tag{6.5a}$$

$$C_{D,p} = \oint_{\partial\Omega} \frac{p\,\boldsymbol{n}\cdot\boldsymbol{t}}{0.5\rho_\infty U_\infty^2 L_\infty}\,\mathrm{d}S \tag{6.5b}$$

$$C_{D,\mu} = \oint_{\partial\Omega} \frac{\boldsymbol{\tau}_w\cdot\boldsymbol{t}}{0.5\rho_\infty U_\infty L_\infty}\,\mathrm{d}S \tag{6.5c}$$

where $\boldsymbol{t} = \big(\cos(\alpha), \sin(\alpha)\big)$ and $\boldsymbol{\tau}_w$ is the wall shear stress vector which is obtained from contracting by the viscous stress tensor $\underline{\underline{\tau}}$ with the surface normal $\hat{\boldsymbol{n}} = -\boldsymbol{n}$:

$$\boldsymbol{\tau}_w = \underline{\underline{\tau}}\,\hat{\boldsymbol{n}} = \begin{pmatrix}\tau_{xx} & \tau_{xy}\\ \tau_{yx} & \tau_{yy}\end{pmatrix}\begin{pmatrix}-n_x\\ -n_y\end{pmatrix}. \tag{6.6}$$

We observe excellent agreement with reference data, despite the relatively coarse grid resolution. Plots of the oscillating (due to the offstream of vortices from the trailing edge, cf. Fig. 9) lift and drag coefficients are shown in Fig. 10.

### *6.5. Sound Generation of a co-rotating Vortex Pair*

For a demonstration of the time integration scheme in a multi-physics application, we apply the P-ERK method in the simulation of the sound generation of a co-rotating vortex pair, which creates a spiral noise pattern [88, 89]. We use a direct-hybrid simulation approach [90], where the flow field is predicted by a compressible Euler simulation, from which acoustic source terms are extracted. These source terms are then fed into a second, coupled simulation using the acoustic perturbation equations to predict the noise propagation. Both the flow and acoustic simulations are performed simultaneously on a shared Cartesian quadtree mesh and exchange source term data at each time step.

| Source | $\overline{C}_L$ | $\overline{C}_D$ |
|---|---|---|
| P-ERK $p = 4, k = 3$ | 0.3827 | 0.04995 |
| P-ERK $p = 2, k = 3$ [1] | 0.3841 | 0.04990 |
| P-ERK $p = 3, k = 3$ [2] | 0.3848 | 0.04910 |
| Uranga et al. [86] | 0.3755 | 0.04978 |
| López-Morales et al. [87] | 0.3719 | 0.04940 |

Table 8: Time-averaged lift and drag coefficients $\overline{C}_L$ and $\overline{C}_D$ for the SD7003 airfoil over the $[30t_c, 35t_c]$ interval.

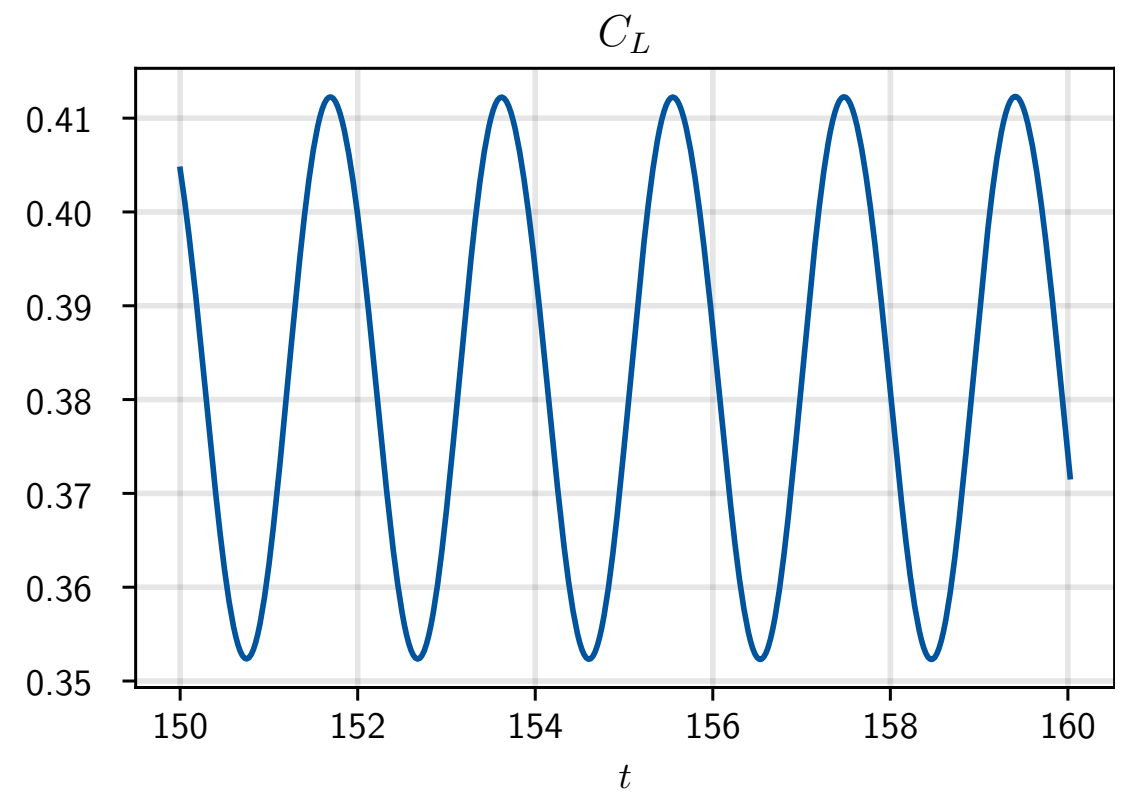

(a) Lift coefficient $C_L = C_{L,p}$, cf. (6.4).

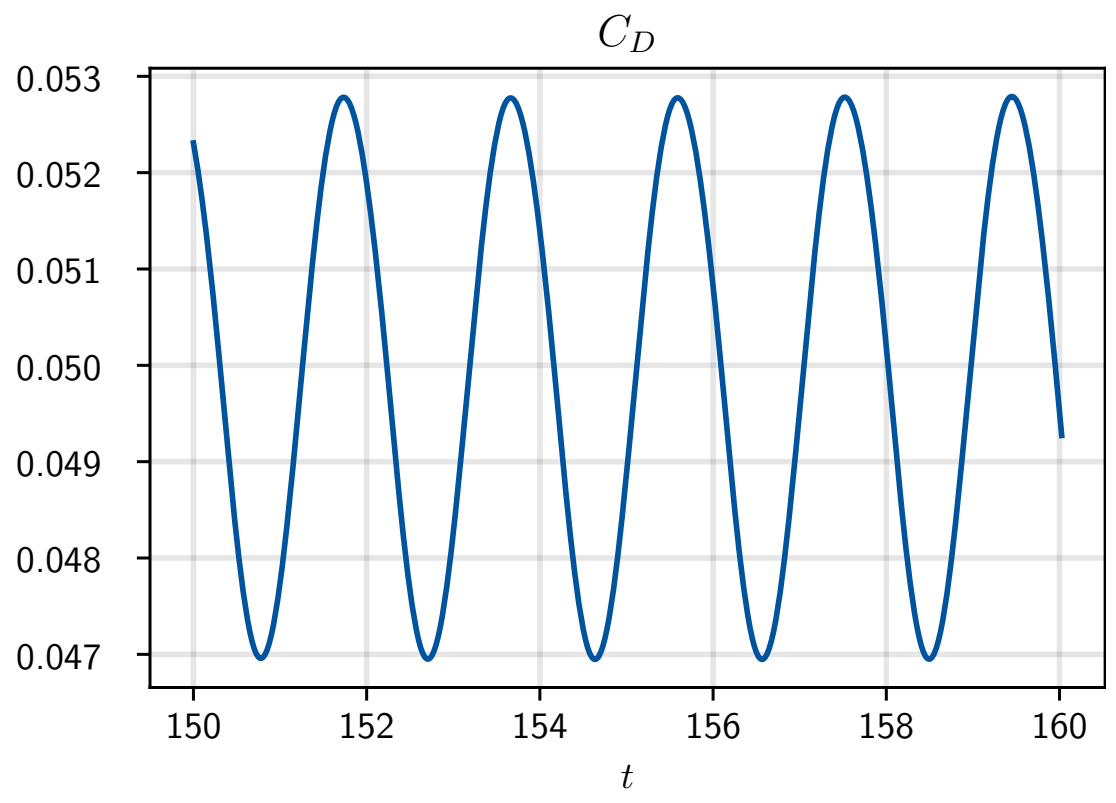

(b) Drag coefficient $C_D = C_{D,p} + C_{D,\mu}$, cf. (6.5).

Figure 10: Unsteady lift (10a) and drag coefficient (10b) for the SD7003 airfoil over the $[30t_c, 32t_c]$ time interval.

The two vortices of equal strength $\Gamma$ and core radius $r_c$ are separated by a distance of $2r_0$ and rotate around each other in a circle, as illustrated in Fig. 11a. The rotational Mach number is $M_r = \frac{1}{9}$, the Reynolds number $\text{Re} = 1.14 \cdot 10^5$, and the core radius is $r_c/r_0 = \frac{2}{9}$. Under the assumption of inviscid and incompressible flow, an analytical solution for this problem can be derived [91]. The initial velocity field is set to the analytical solution, while the density is assumed constant. The hydrodynamic pressure is determined by the Bernoulli equation. All perturbed quantities in the acoustics simulation are initially set to zero. For an in-depth description of the numerical setup, please refer to [90, 91].

The computational domain is a square with side length $270r_0$, with circular refinement patches towards the center, see Fig. 11b. In the center with the hightest mesh resolution, the acoustic source terms are calculated from the velocity field of the flow solver and subsequently used in the acoustics solver. Both for the flow solver and the acoustics solver, we use local solution polynomials of degree $k = 3$ and the Rusanov/local Lax-Friedrichs approximate Riemann solver. The domain is initially discretized using a $128^2$ mesh with four concentric spherical refinements, cf. Fig. 11b. For the flow solver, subsonic outflow conditions are prescribed. In addition, a sponge layer with thickness $20r_0$ is applied at the domain boundaries, which dampens density and pressure fluctuations towards the freestream values. In the CAA simulation, zero values are used as the boundary state for the perturbed quantities and a sponge layer with thickness $35r_0$ along the domain boundary dampens the perturbed pressure towards zero. Acoustic source terms are only calculated on the central, circular refinement patch defined by $r/r0 < 5$, see Fig. 11b.

We optimize two families of P-ERK schemes for the sound generation by a co-rotating vortex pair. First, we construct a four member family with $E = \{5, 6, 8, 13\}$ that is optimized for the Euler-only semidiscretization. Additionally, we construct a four member family with $E = \{5, 6, 9, 14\}$ that is optimized for the Euler-acoustic semidiscretization incorporating the source term from the Euler equations. As for the previous example, we compare the optimized

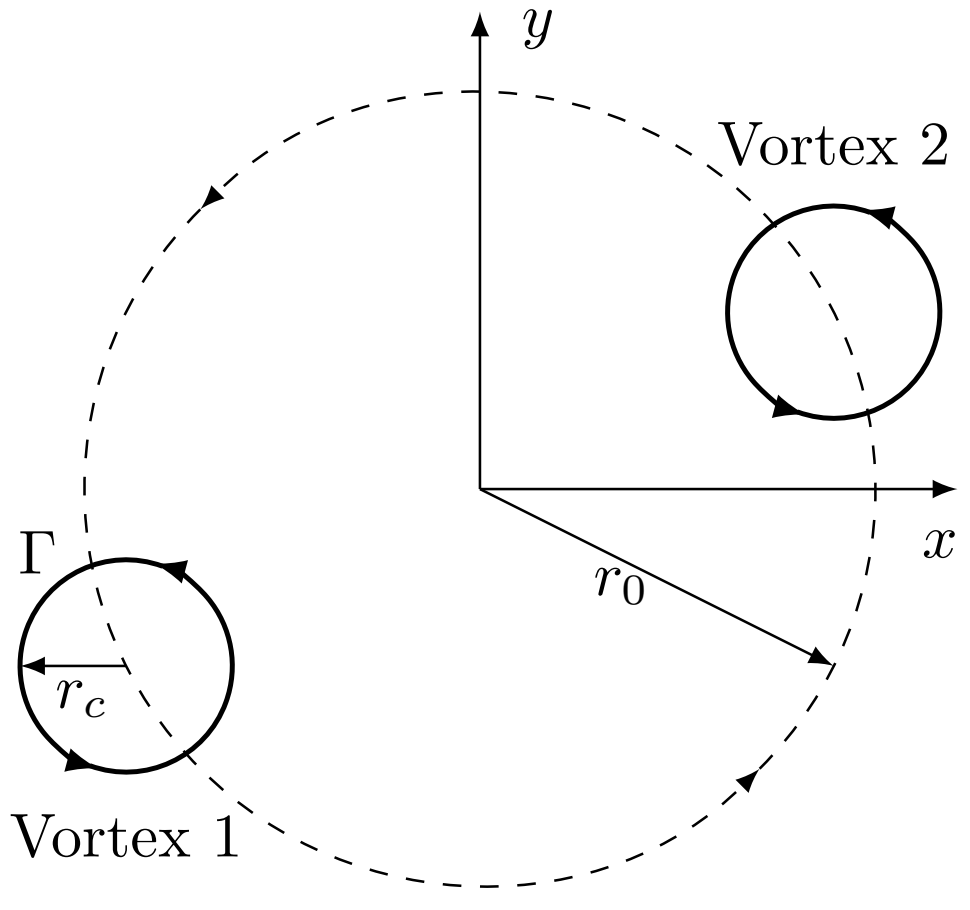

(a) Flow setup for the spinning vortices located on a circle with radius $r_0$, core radius $r_c$, and strength $\Gamma$.

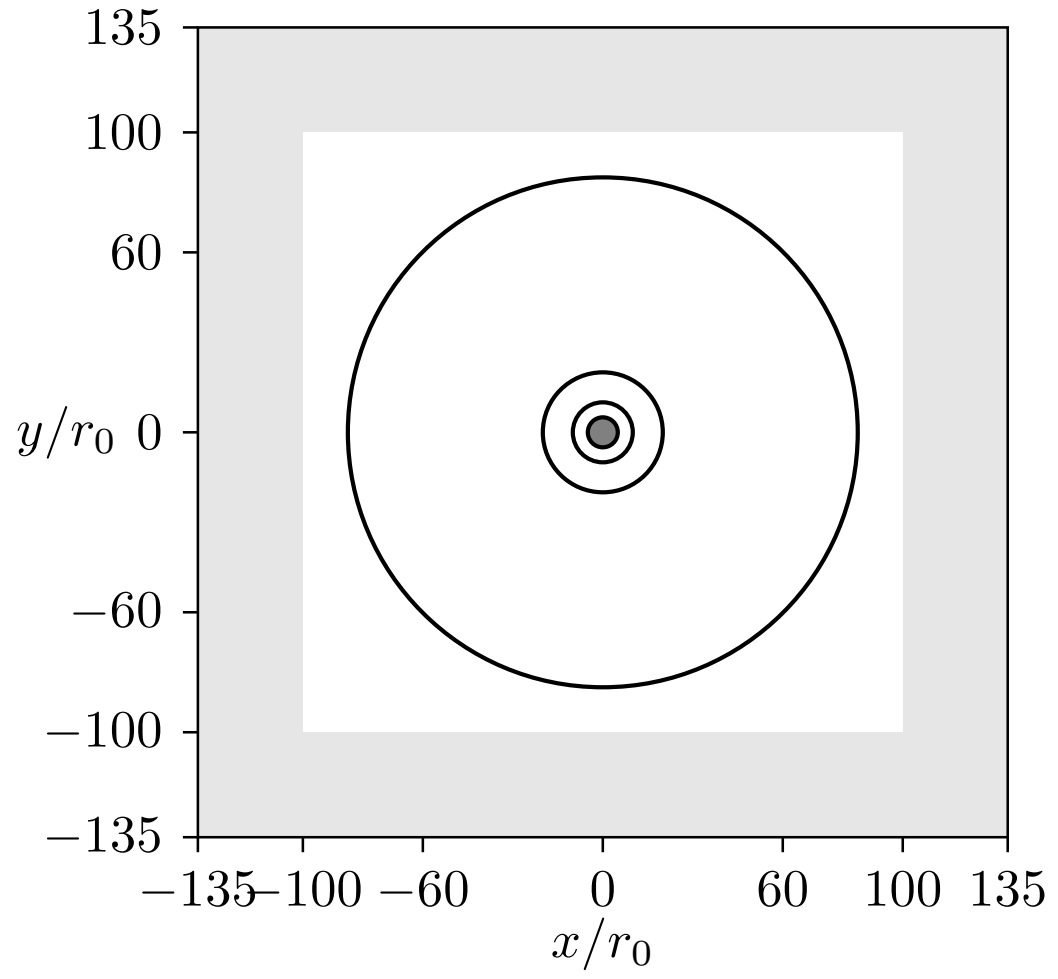

(b) Domain with refinement level jumps (circles) towards the center, where the solvers are coupled via the source terms (dark shaded region). The sponge-layer is indicated by the light shaded region.

Figure 11: Numerical setup for the coupled flow-acoustics simulation of two co-rotating vortices. Figures are based on [90].

| Method | $\tau/\tau^\star$ | $N_{\text{RHS}}/N_{\text{RHS}}^\star$ |
|---|---|---|
| P-ERK$_{4;\{5,6,8,13\}}$, P-ERK$_{4;\{5,6,9,14\}}$ | 1.0 | 1.0 |
| P-ERK$_{4;13}$, P-ERK$_{4;14}$ | 1.61 | 1.87 |
| NDB$_{4;14}$ [77] | 1.82 | 2.20 |
| TD$_{4;8}$ [79] | 2.65 | 3.02 |
| CFR$_{4;6}$ [54] | 3.86 | 4.26 |
| RK$_{4;4}$ | 4.77 | 4.52 |

Table 9: Runtimes and number of scalar RHS evaluations (6.3) of different fourth-order integrators compared to the optimized $p = 4, E = \{5, 6, 8, 13\}$, $p = 4, E = \{5, 6, 9, 14\}$ integrators for the sound generation by the vortex pair.

P-ERK schemes to the standalone P-ERK$_{4;13}$, P-ERK$_{4;14}$ schemes, the canonical RK$_{4;4}$ scheme, and the well-performing NDB$_{4;14}$ scheme. Furthermore, we also run this example with the six-stage, low-dissipation, low-dispersion method CFR$_{4;6}$ [54] optimized for computational acoustics and the eight-stage, low-storage, low-dissipation, low-dispersion scheme TD$_{4;8}$ [79] optimized for DG methods applied to wave propagation problems.

For this example, we observe the most pronounced speedup for the multirate P-ERK schemes compared to the other methods. The standalone P-ERK schemes require 87% more scalar RHS evaluations which manifests in 61% higher runtime. Among the off-the-shelf methods from `OrdinaryDiffEq.jl`, the NDB$_{4;14}$ scheme performs best, being 82% slower than the multirate P-ERK schemes and requiring more than twice the scalar RHS evaluations. As for the first example discussed in this section, the canonical RK$_{4;4}$ is again disproportionally slow relative to the ratio of scalar RHS evaluations. The resulting, well-known spiral noise pattern [88, 89] of this setup is shown in Fig. 12.

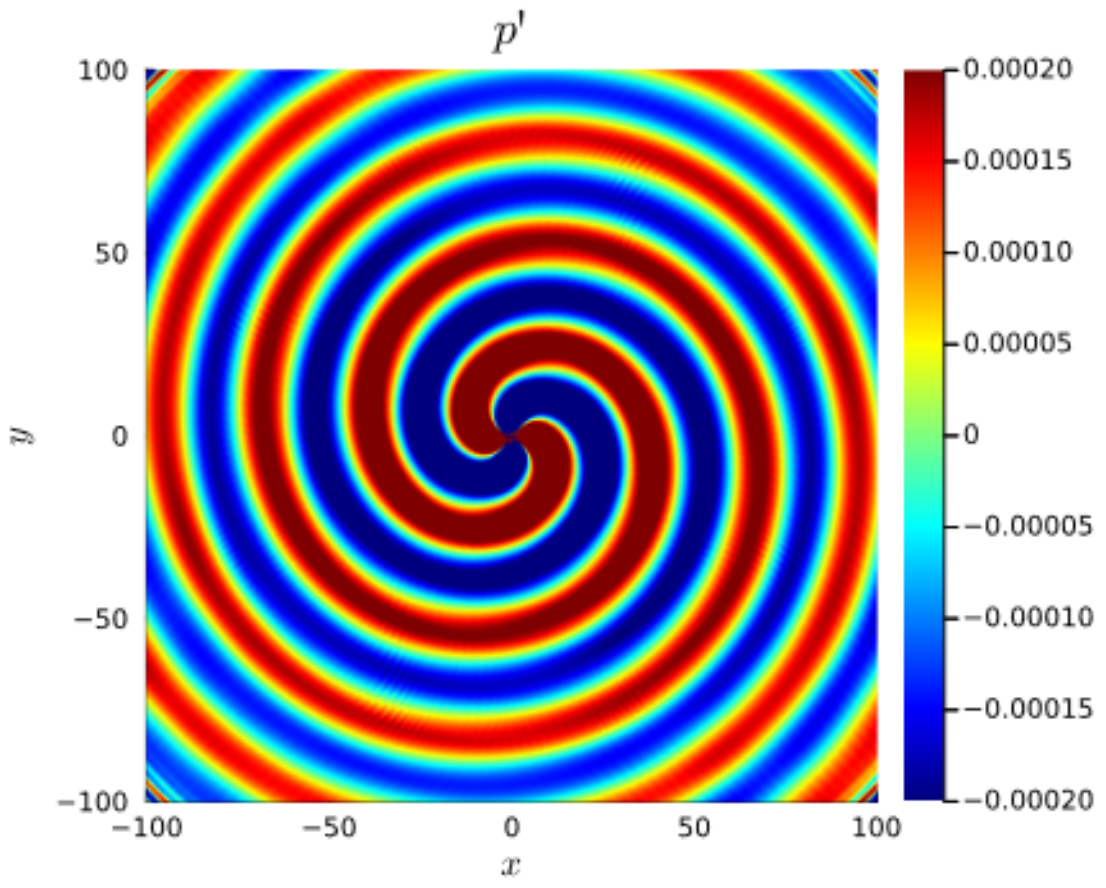

Figure 12: Pressure perturbation $p'$ at $t = 5t_a$ for the sound generation by the vortex pair in the non-damped region, cf. Fig. 11b.

## 7. Conclusions

In this paper, we propose a novel family of partitioned explicit Runge-Kutta schemes extending the previous works on second- [1] and third-order [2] P-ERK schemes to fourth order. The P-ERK schemes offer a seamless multirate time integration mechanism for partitioned problems with different characteristic speeds. By honoring the order-condition that couples Butcher arrays of different schemes, we motivate a certain structure of the Butcher arrays which fulfills this condition by construction. Building upon this particular structure, we derive a new parametrization of the corresponding stability polynomials which is required for an efficient optimization of the schemes' domain of stability using the algorithm developed in [51]. We thoroughly validated the herein derived schemes in terms of linear stability, pointwise errors, conservation of linear invariants, and order of convergence.

We demonstrated the efficiency of the optimized schemes for a variety of applications with spatially varying characteristic speeds. These include the linearized Euler equations with non-uniform mean speed of sound, where stiffness is induced due to a multiscale problem setup. Building upon a previous work [43], we also studied the dynamic re-partitioning by means of the isentropic vortex advection on an adaptively refined mesh. Considering a completely different flow regime, we applied the fourth-order P-ERK schemes to the laminar flow around the SD7003 airfoil. Finally, we studied the sound generation by a co-rotating vortex pair in a hybrid simulation setup which benefits heavily from the optimized schemes.

For all these applications, the multirate P-ERK schemes outperform the standalone optimized schemes and a range of state-of-the art methods available in `OrdinaryDiffEq.jl`. For these, the RHS evaluations are cut down by factors two to three, reducing the observed runtimes by 40% to 280%. Compared to the canonical fourth-order Runge-Kutta method, the multirate P-ERK schemes may reduce runtime and RHS evaluations by more than a factor of four.

In future work, we will explore P-ERK schemes with error-based timestep control by means of third- and fourth-order methods, thereby extending [44].

## Data Availability

All data generated or analyzed during this study are included in this published article and its supplementary information files.

## Code Availability & Reproducibility

We will provide a reproducibility repository in a later version of this article. The current progress can be found in the GitHub repository

`https://github.com/DanielDoehring/Trixi.jl/tree/AdaptiveTimeIntegration_p4`.

**Acknowledgments**

Funding by German Research Foundation (DFG) under Research Unit FOR5409: "Structure-Preserving Numerical Methods for Bulk- and Interface-Coupling of Heterogeneous Models (SNuBIC) " (grant #463312734).

**Declaration of competing interest**

The authors declare the following financial interests/personal relationships which may be considered as potential competing interests: Daniel Doehring's financial support was provided by German Research Foundation.

**CRediT authorship contribution statement**

**Daniel Doehring**: Formal analysis, Investigation, Methodology, Conceptualization, Software, Writing - original draft.
**Lars Christmann**: Software, Writing - review & editing.
**Michael Schlottke-Lakemper**: Software, Conceptualization, Funding acquisition, Writing - review & editing.
**Gregor J. Gassner**: Conceptualization, Funding acquisition, Writing - review & editing.
**Manuel Torrilhon**: Conceptualization, Funding acquisition, Supervision, Writing - review & editing.

## Appendix A. Order Conditions for P-ERK

Here, we give all eight order constraints following from (2.3), (2.4) for the P-ERK weight vector $\boldsymbol{b}$ (3.6) and the archetype Butcher matrix $A^{(r)}$ (3.7):

$$p = 1: \quad 1 \overset{!}{=} b_{S-1} + b_S \tag{A.1a}$$

$$p = 2: \quad \frac{1}{2} \overset{!}{=} b_{S-1}\, c_{S-1} + b_S\, c_S \tag{A.1b}$$

$$p = 3: \quad \frac{1}{3} \overset{!}{=} b_{S-1}\, c_{S-1}^2 + b_S\, c_S^2 \tag{A.1c}$$

$$\frac{1}{6} \overset{!}{=} b_{S-1}\, a_{S-1,S-2}^{(r)}\, c_{S-1} + b_S\, a_{S,S-1}^{(r)}\, c_S \qquad r = 1, \dots, R \tag{A.1d}$$

$$p = 4: \quad \frac{1}{4} \overset{!}{=} b_{S-1}\, c_{S-1}^3 + b_S\, c_S^3 \tag{A.1e}$$

$$\frac{1}{8} \overset{!}{=} b_{S-1}\, c_{S-1}\, a_{S-1,S-2}^{(r)}\, c_{S-2} + b_S\, c_S\, a_{S,S-1}^{(r)}\, c_{S-1} \qquad r = 1, \dots, R \tag{A.1f}$$

$$\frac{1}{12} \overset{!}{=} b_{S-1}\, a_{S-1,S-2}^{(r)}\, c_{S-2}^2 + b_S\, a_{S,S-1}^{(r)}\, c_{S-1}^2 \qquad r = 1, \dots, R \tag{A.1g}$$

$$\frac{1}{24} \overset{!}{=} b_{S-1}\, a_{S-1,S-2}^{(r_1)}\, a_{S-2,S-3}^{(r_2)}\, c_{S-3} + b_S\, a_{S,S-1}^{(r_1)}\, a_{S-1,S-2}^{(r_2)}\, c_{S-2} \qquad r_1, r_2 = 1, \dots, R \tag{A.1h}$$

## Appendix B. Shared Stages: Influence on Nonlinear Stability

In this appendix we seek to demonstrate the benefits of shared stages, i.e., the combination of a true partitioned and classic Runge-Kutta method for the nonlinear stability properties of the P-ERK schemes. The analysis here is inspired by [65] where we conduct a thorough investigation of the nonlinear stability properties of the P-ERK schemes.

We consider the 1D advection equation

$$u_t + u_x = 0, \tag{B.1a}$$

$$u_0(t_0 = 0, x) = 1 + \frac{1}{2}\sin(\pi x) \tag{B.1b}$$

on domain $\Omega = (-1, 1)$ equipped with periodic boundary conditions. The cell averages

$$U_i(t) := \frac{1}{\Delta x_i}\int_{x_{i-1/2}}^{x_{i+1/2}} u(t,x)\,\mathrm{d}x, \quad i = 1,\dots,N \tag{B.2}$$

are evolved using the Godunov/upwind flux on a non-uniform mesh. Here, we consider the second-order P-ERK method with $E = \{8, 16\}$ stage evaluations. The domain is refined in the center $(-0.5, 0.5)$ such that

$$\Delta x_i = \begin{cases} \frac{1}{32} & i = 1,\dots,16 \\ \frac{1}{64} & i = 17,\dots,80 \\ \frac{1}{32} & i = 81,\dots 96 \end{cases}. \tag{B.3}$$

The $E = 16$ stage-evaluation method is used in the center of the domain, i.e., for $U_i, i = 17,\dots,80$, and the $E = 8$ stage-evaluation method is used in the remainder of the domain.

We construct two sets of stability polynomials and corresponding P-ERK methods for this case. First, a standard P-ERK scheme cf. (3.5) from (apart from second-order accurate) unconstrained stability polynomials. Second, a specialized P-ERK scheme where the last six stages are shared which is achieved by restricting the first eight monomial coefficients $\alpha_i, i = 1,\dots 8$ to be equal between the $E = 8$ and $E = 16$ stability polynomial. For the second-order accurate stability polynomials this implies that six in principle free coefficients are now shared between methods. This reduces the optimal timestep of this method to about 64.3% of the optimal, unconstrained timestep. The structure of the Butcher tableaus for the non-standard method is similar to (3.10) and in this case given by

$$\begin{array}{c|c|cc} & & \multicolumn{2}{c}{r} \\ \hline i & \boldsymbol{c} & 1 & 2 \\ \cline{3-4} & & A^{(1)} & A^{(2)} \\ \hline & & \multicolumn{2}{c}{\boldsymbol{b}^T} \end{array} = \begin{array}{c|c|cc} & & \multicolumn{2}{c}{r} \\ \hline i & \boldsymbol{c} & 1 & 2 \\ \hline 1 & 0 & \multicolumn{2}{c}{\multirow{2}{*}{$A_{1:2}$}} \\ 2 & c_2 & & \\ \hdashline 3 & c_3 & & \\ \vdots & \vdots & A^{(1)}_{3:10} & A^{(2)}_{3:10} \\ 10 & c_{10} & & \\ \hdashline 11 & c_{11} & & \\ \vdots & \vdots & \multicolumn{2}{c}{A_{11:16}} \\ 16 & c_{16} & & \\ \hline & & \multicolumn{2}{c}{\boldsymbol{b}^T} \end{array} \tag{B.4}$$

We present $\boldsymbol{U}$ after one timestep of this size for both standard and shared-stage P-ERK scheme in Fig. B.13. Clearly, the shared-stage P-ERK scheme leads to significantly less severe oscillations compared to the standard P-ERK scheme at the interface between the two methods. This is noteworthy, as the standard scheme allows for a significantly larger linearly stable timestep.

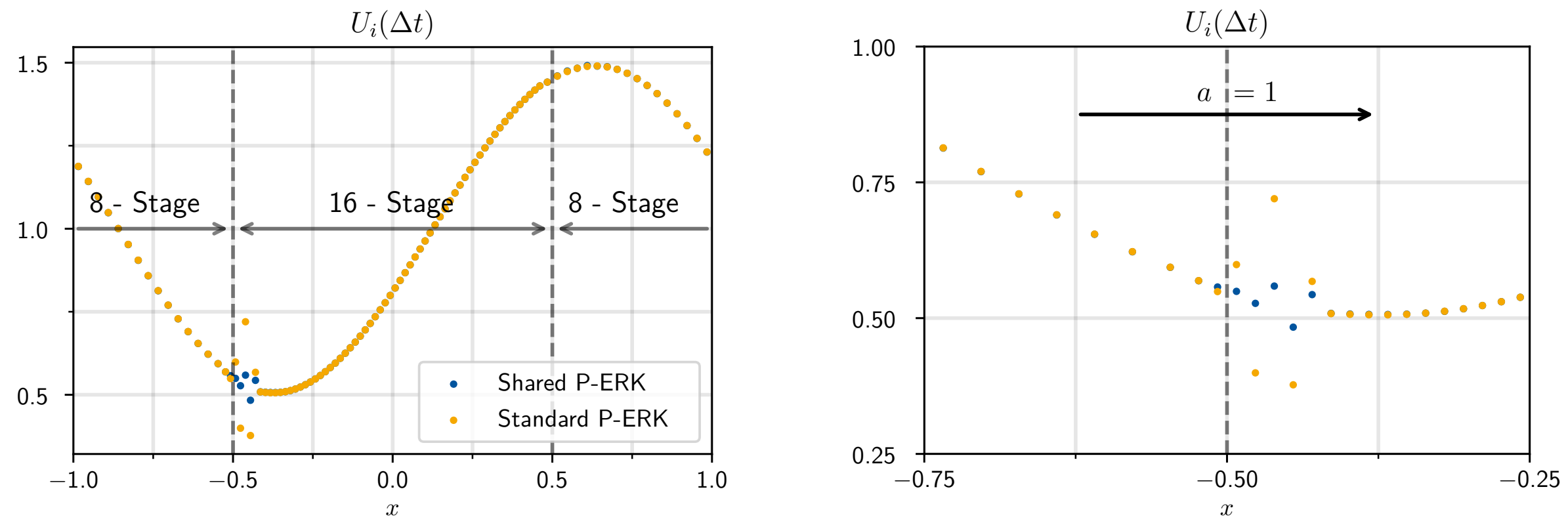

(a) Entire domain $\Omega = (-1, 1)$.

(b) Close view of the scheme interface at $x = -0.5$.

Figure B.13: Cell averages $\boldsymbol{U}$ after one timestep of size $\Delta t \approx 0.14067$ for two second-order P-ERK scheme with $E = \{8, 16\}$ stage evaluations.

## Appendix C. Second-Order Cone Formulation of the Linear Stability Constraint

To derive the second-order cone form of the stability constraint (4.1), we start by noting that we can write (4.1) for stability polynomials $P_{4;E}(z;\boldsymbol{\gamma})$ of form (4.7) as

$$\left|\boldsymbol{Z}_m \boldsymbol{\gamma} + b_m\right| \leq 1 \quad m = 1, \dots, M \tag{C.1}$$

with $\boldsymbol{\gamma} \in \mathbb{R}^{E-5}$ and

$$\boldsymbol{Z}_m^T \coloneqq \begin{pmatrix} c_{E-4}\left(\frac{k_2}{c_{E-3}} z^5 + \frac{k_1}{c_{E-3}} z^6\right) \\ c_{E-5}\left(\frac{k_2}{c_{E-3}} z^6 + \frac{k_1}{c_{E-3}} z^7\right) \\ \vdots \\ c_2\left(\frac{k_2}{c_{E-3}} z^{E-1} + \frac{k_1}{c_{E-3}} z^E\right) \end{pmatrix} \in \mathbb{C}^{E-5}, \quad \boldsymbol{Z}_m \in \mathbb{C}^{1\times(E-5)} \tag{C.2a}$$

$$b_m \coloneqq \sum_{j=0}^{p} \frac{z_m^j}{j!} + k_1 z^5 \in \mathbb{C}\,. \tag{C.2b}$$

By splitting the $\boldsymbol{Z}_m, b_m$ into real and imaginary parts, the left-hand-side of (C.1) can be rewritten as

$$\left|\boldsymbol{Z}_m \boldsymbol{\gamma} + b_m\right| = \left|\mathrm{Re}(\boldsymbol{Z}_m)\boldsymbol{\gamma} + \mathrm{Re}(b_m) + i\Big(\mathrm{Im}(\boldsymbol{Z}_m)\boldsymbol{\gamma} + \mathrm{Im}(b_m)\Big)\right| \tag{C.3}$$

$$= \sqrt{\Big(\mathrm{Re}(\boldsymbol{Z}_m)\boldsymbol{\gamma} + \mathrm{Re}(b_m)\Big)^2 + \Big(\mathrm{Im}(Z_m)\boldsymbol{\gamma} + \mathrm{Im}(b_m)\Big)^2} \tag{C.4}$$

$$= \left\| \begin{pmatrix} \mathrm{Re}(\boldsymbol{Z}_m) \\ \mathrm{Im}(\boldsymbol{Z}_m) \end{pmatrix} \boldsymbol{\gamma} + \begin{pmatrix} \mathrm{Re}(b_m) \\ \mathrm{Im}(b_m) \end{pmatrix} \right\|_2 \tag{C.5}$$

$$=: \left\|\widetilde{Z}_m \boldsymbol{\gamma} + \widetilde{\boldsymbol{b}}_m\right\|_2 \tag{C.6}$$

with $\widetilde{Z}_m \in \mathbb{R}^{2\times(E-5)}, \widetilde{\boldsymbol{b}}_m \in \mathbb{R}^2$. Through these transformations, the stability constraint (4.1) can be reformulated as a *Second-Order Cone Constraint* [92, 93] of the form

$$\left\|\widetilde{Z}_m \boldsymbol{\gamma} + \widetilde{\boldsymbol{b}}_m\right\|_2 \leq 1 \quad m = 1, \dots, M. \tag{C.7}$$

Further, the optimization problem (4.2) can be equivalently written as

$$\min_{\widetilde{\gamma} \in \mathbb{R}^{S-4}} \boldsymbol{d}^T \widetilde{\boldsymbol{\gamma}} \text{ such that } \left\| \begin{pmatrix} 0 \\ \begin{pmatrix} \boldsymbol{0} & \widetilde{Z}_m \end{pmatrix} \end{pmatrix} \widetilde{\boldsymbol{\gamma}} + \begin{pmatrix} 0 \\ \widetilde{\boldsymbol{b}}_m \end{pmatrix} \right\|_2 \leq \boldsymbol{d}^T \widetilde{\boldsymbol{\gamma}} \quad m = 1, \dots, M \tag{C.8}$$

where $\widetilde{\boldsymbol{\gamma}} \in \mathbb{R}^{S-4}$, since an additional variable indicating the value of the stability polynomial at most violating eigenvalue is added to $\boldsymbol{\gamma}$. The objective and simultaneous constraint vector $\boldsymbol{d} = \begin{pmatrix} 1 & 0 & \dots & 0 \end{pmatrix} \in \mathbb{R}^{E-4}$ ensures that for a given $\Delta t$ the $\boldsymbol{\gamma}$ are chosen such that the constraint violation $\widetilde{\gamma}_1$ is minimized. (C.8) is a *second-order Cone Problem* [93, 92] for which specifically tailored solvers are available.

Here, `EiCOS` [94], a `c++` adaptation of `ECOS` [95] is used. The reason to not use, e.g., the `Matlab`-based `CVX` interface [96, 97] as done in [51, 1, 2], is that `Eigen3` [98] is written in a type-generic (`c++ templates`) fashion and thus allows usage of data types with precision beyond standard double precision accuracy. In particular, we implemented derivates of the `EiCOS` solver which allows for the use of the `c++` built-in `long double` floating point datatype as well as arbitrary precision data types from the `boost` multiprecision library [99]. The optimization problem (C.8) can be implemented in about 200 lines of code, which is made publicly available on GitHub [100].

## References

[1] B. C. Vermeire, Paired explicit Runge-Kutta schemes for stiff systems of equations, Journal of Computational Physics 393 (2019) 465–483. doi:`10.1016/j.jcp.2019.05.014`.

[2] S. H. Nasab, B. C. Vermeire, Third-order paired explicit Runge-Kutta schemes for stiff systems of equations, Journal of Computational Physics 468 (2022) 111470. doi:`10.1016/j.jcp.2022.111470`.

[3] M. Günther, A. Kvaernø, P. Rentrop, Multirate partitioned Runge-Kutta methods, BIT Numerical Mathematics 41 (2001) 504–514. doi:`10.1023/A:1021967112503`.

[4] E. M. Constantinescu, A. Sandu, Multirate timestepping methods for hyperbolic conservation laws, Journal of Scientific Computing 33 (2007) 239–278. doi:`10.1007/s10915-007-9151-y`.

[5] M. Günther, A. Sandu, Multirate generalized additive Runge-Kutta methods, Numerische Mathematik 133 (2016) 497–524. doi:`10.1007/s00211-015-0756-z`.

[6] S. Osher, R. Sanders, Numerical approximations to nonlinear conservation laws with locally varying time and space grids, Mathematics of computation 41 (1983) 321–336. doi:`10.1090/S0025-5718-1983-0717689-8`.

[7] H.-Z. Tang, G. Warnecke, High resolution schemes for conservation laws and convection-diffusion equations with varying time and space grids, Journal of computational mathematics (2006) 121–140. URL: `https://www.jstor.org/stable/43694072`.

[8] M. J. Berger, J. Oliger, Adaptive mesh refinement for hyperbolic partial differential equations, Journal of Computational Physics 53 (1984) 484–512. doi:`10.1016/0021-9991(84)90073-1`.

[9] M. J. Grote, M. Mehlin, T. Mitkova, Runge–Kutta-based explicit local time-stepping methods for wave propagation, SIAM Journal on Scientific Computing 37 (2015) A747–A775. doi:`10.1137/140958293`.

[10] F. Lörcher, G. Gassner, C.-D. Munz, A discontinuous Galerkin scheme based on a space–time expansion. I. inviscid compressible flow in one space dimension, Journal of Scientific Computing 32 (2007) 175–199. doi:10.1007/s10915-007-9128-x.

[11] L. Krivodonova, An efficient local time-stepping scheme for solution of nonlinear conservation laws, Journal of Computational Physics 229 (2010) 8537–8551. doi:10.1016/j.jcp.2010.07.037.

[12] M. Dumbser, M. Käser, E. F. Toro, An arbitrary high-order discontinuous Galerkin method for elastic waves on unstructured meshes–V. local time stepping and p-adaptivity, Geophysical Journal International 171 (2007) 695–717. doi:10.1111/j.1365-246X.2007.03427.x.

[13] C. Dawson, R. Kirby, High resolution schemes for conservation laws with locally varying time steps, SIAM Journal on Scientific Computing 22 (2001) 2256–2281. doi:10.1137/S1064827500367737.

[14] V. Kulka, P. Jenny, Temporally adaptive conservative scheme for unsteady compressible flow, Journal of Computational Physics 455 (2022) 110918. doi:10.1016/j.jcp.2021.110918.

[15] B. Cockburn, C.-W. Shu, Runge–Kutta discontinuous Galerkin methods for convection-dominated problems, Journal of Scientific Computing 16 (2001) 173–261. doi:10.1023/A:1012873910884.

[16] U. M. Ascher, S. J. Ruuth, B. T. Wetton, Implicit-explicit methods for time-dependent partial differential equations, SIAM Journal on Numerical Analysis 32 (1995) 797–823. doi:10.1137/0732037.

[17] E. Hairer, Order conditions for numerical methods for partitioned ordinary differential equations, Numerische Mathematik 36 (1981) 431–445. doi:10.1007/BF01395956.

[18] W. Hundsdorfer, D. I. Ketcheson, I. Savostianov, Error analysis of explicit partitioned Runge-Kutta schemes for conservation laws, Journal of Scientific Computing 63 (2015) 633–653. doi:10.1007/s10915-014-9906-1.

[19] I. Higueras, Strong stability for additive Runge–Kutta methods, SIAM Journal on Numerical Analysis 44 (2006) 1735–1758. doi:10.1137/040612968.

[20] I. Higueras, Characterizing strong stability preserving additive Runge-Kutta methods, Journal of Scientific Computing 39 (2009) 115–128. doi:10.1007/s10915-008-9252-2.

[21] W. Hundsdorfer, A. Mozartova, V. Savcenco, Analysis of explicit multirate and partitioned Runge-Kutta schemes for conservation laws, Technical Report, Centrum voor Wiskunde en Informatica, 2007. URL: https://research.tue.nl/en/publications/analysis-of-explicit-multirate-and-partitioned-{R}unge-{K}utta-scheme, cWI report. MAS-E; Vol. 0715.

[22] D. I. Ketcheson, C. B. MacDonald, S. J. Ruuth, Spatially partitioned embedded Runge–Kutta methods, SIAM Journal on Numerical Analysis 51 (2013) 2887–2910. doi:10.1137/130906258.

[23] E. Hairer, G. Wanner, S. P. Nørsett, Solving Ordinary Differential Equations I: Nonstiff Problems, Springer Series in Computational Mathematics, 2 ed., Springer Berlin, Heidelberg, 1993. doi:10.1007/978-3-540-78862-1.

[24] E. Hofer, A partially implicit method for large stiff systems of odes with only few equations introducing small time-constants, SIAM Journal on Numerical Analysis 13 (1976) 645–663. doi:10.1137/0713054.

[25] E. Griepentrog, Gemischte Runge-Kutta-verfahren für steife systeme, Seminarbereicht Sekt. Math (1978) 19–29.

[26] L. Jay, Symplectic partitioned Runge–Kutta methods for constrained Hamiltonian systems, SIAM Journal on Numerical Analysis 33 (1996) 368–387. doi:10.1137/0733019. arXiv:10.1137/0733019.

[27] L. Abia, J. Sanz-Serna, Partitioned Runge-Kutta methods for separable Hamiltonian problems, mathematics of Computation 60 (1993) 617–634. doi:10.1090/S0025-5718-1993-1181328-1.

[28] C. A. Kennedy, M. H. Carpenter, Additive Runge–Kutta schemes for convection–diffusion–reaction equations, Applied Numerical Mathematics 44 (2003) 139–181. doi:10.1016/S0168-9274(02)00138-1.

[29] G. J. Cooper, A. Sayfy, Additive Runge-Kutta methods for stiff ordinary differential equations, Mathematics of Computation 40 (1983) 207–218. doi:10.1090/S0025-5718-1983-0679441-1.

[30] A. Sandu, M. Günther, A generalized-structure approach to additive Runge–Kutta methods, SIAM Journal on Numerical Analysis 53 (2015) 17–42. doi:10.1137/130943224.

[31] P. Rentrop, Partitioned Runge-Kutta methods with stiffness detection and stepsize control, Numerische Mathematik 47 (1985) 545–564. doi:10.1007/BF01389456.

[32] J. Bruder, K. Strehmel, R. Weiner, Partitioned adaptive Runge-Kutta methods for the solution of nonstiff and stiff systems, Numerische Mathematik 52 (1988) 621–638. doi:10.1007/BF01395815.

[33] P. Albrecht, A new theoretical approach to Runge–Kutta methods, SIAM Journal on Numerical Analysis 24 (1987) 391–406. doi:10.1137/0724030.

[34] Z. Jackiewicz, R. Vermiglio, Order conditions for partitioned Runge-Kutta methods, Applications of Mathematics 45 (2000) 301–316. URL: http://dml.cz/dmlcz/134441.

[35] J. D. Lambert, Numerical methods for ordinary differential systems, Wiley New York, 1992. URL: https://www.wiley.com/en-us/Numerical+Methods+for+Ordinary+Differential+Systems%3A+The+Initial+Value+Problem-p-9780471929901.

[36] G. Wanner, E. Hairer, Solving ordinary differential equations II: Stiff and Differential-Algebraic Problems, volume 375 of *Springer Series in Computational Mathematics*, 2 ed., Springer Berlin, Heidelberg, 1996. doi:10.1007/978-3-642-05221-7.

[37] R. I. McLachlan, Y. Sun, P. Tse, Linear stability of partitioned Runge–Kutta methods, SIAM Journal on Numerical Analysis 49 (2011) 232–263. doi:10.1137/100787234.

[38] J. P. Hespanha, Linear systems theory, 2 ed., Princeton University press, 2018. URL: https://press.princeton.edu/books/hardcover/9780691179575/linear-systems-theory.

[39] C.-W. Shu, S. Osher, Efficient implementation of essentially non-oscillatory shock-capturing schemes, Journal of Computational Physics 77 (1988) 439–471. doi:10.1016/0021-9991(88)90177-5.

[40] S. Gottlieb, D. Ketcheson, C.-W. Shu, Strong stability preserving Runge-Kutta and multistep time discretizations, World Scientific, 2011. doi:10.1142/9789814289276_0001.

[41] J. F. B. M. Kraaijevanger, Contractivity of Runge-Kutta methods, BIT Numerical Mathematics 31 (1991) 482–528. doi:10.1007/BF01933264.

[42] E. J. Kubatko, B. A. Yeager, D. I. Ketcheson, Optimal strong-stability-preserving Runge–Kutta time discretizations for discontinuous Galerkin methods, Journal of Scientific Computing 60 (2014) 313–344. doi:10.1007/s10915-013-9796-7.

[43] D. Doehring, M. Schlottke-Lakemper, G. J. Gassner, M. Torrilhon, Multirate time-integration based on dynamic ode partitioning through adaptively refined meshes for compressible fluid dynamics, Journal of Computational Physics 514 (2024) 113223. doi:10.1016/j.jcp.2024.113223.

[44] B. C. Vermeire, Embedded paired explicit Runge-Kutta schemes, Journal of Computational Physics 487 (2023) 112159. doi:10.1016/j.jcp.2023.112159.

[45] A. Abdulle, Explicit stabilized Runge-Kutta methods, Technical Report, Mathematics Institute of Computational Science and Engineering, School of Basic Sciences - Section of Mathematics EPFL Lausanne, 2011. URL: https://www.epfl.ch/labs/mathicse/wp-content/uploads/2018/10/27.2011_AA.pdf.

[46] P. Van der Houwen, The development of Runge-Kutta methods for partial differential equations, Applied Numerical Mathematics 20 (1996) 261–272. doi:10.1016/0168-9274(95)00109-3.

[47] J. G. Verwer, B. P. Sommeijer, W. Hundsdorfer, RKC time-stepping for advection–diffusion–reaction problems, Journal of Computational Physics 201 (2004) 61–79. doi:10.1016/j.jcp.2004.05.002.

[48] M. Torrilhon, R. Jeltsch, Essentially optimal explicit Runge-Kutta methods with application to hyperbolic–parabolic equations, Numerische Mathematik 106 (2007) 303–334. doi:10.1007/s00211-006-0059-5.

[49] R. Jeltsch, O. Nevanlinna, Largest disk of stability of explicit Runge-Kutta methods, BIT Numerical Mathematics 18 (1978) 500–502. doi:10.1007/BF01932030.

[50] B. Owren, K. Seip, Some stability results for explicit Runge-Kutta methods, BIT Numerical Mathematics 30 (1990) 700–706. doi:10.1007/BF01933217.

[51] D. Ketcheson, A. Ahmadia, Optimal stability polynomials for numerical integration of initial value problems, Communications in Applied Mathematics and Computational Science 7 (2013) 247–271. doi:10.2140/camcos.2012.7.247.

[52] D. I. Ketcheson, L. Lóczi, M. Parsani, Internal error propagation in explicit Runge-Kutta methods, SIAM Journal on Numerical Analysis 52 (2014) 2227–2249. doi:10.1137/130936245.

[53] J. Berland, C. Bogey, C. Bailly, Low-dissipation and low-dispersion fourth-order Runge-Kutta algorithm, Computers & Fluids 35 (2006) 1459–1463. doi:10.1016/j.compfluid.2005.04.003.

[54] M. Calvo, J. Franco, L. Rández, A new minimum storage Runge-Kutta scheme for computational acoustics, Journal of Computational Physics 201 (2004) 1–12. doi:10.1016/j.jcp.2004.05.012.

[55] Z. J. Wang, K. Fidkowski, R. Abgrall, F. Bassi, D. Caraeni, A. Cary, H. Deconinck, R. Hartmann, K. Hillewaert, H. T. Huynh, et al., High-order CFD methods: Current status and perspective, International Journal for Numerical Methods in Fluids 72 (2013) 811–845. doi:10.1002/fld.3767.

[56] K. Black, A conservative spectral element method for the approximation of compressible fluid flow, Kybernetika 35 (1999) 133–146.

[57] D. A. Kopriva, Implementing spectral methods for partial differential equations: Algorithms for scientists and engineers, Scientific Computation (SCIENTCOMP), 1 ed., Springer Science & Business Media, 2009. doi:10.1007/978-90-481-2261-5.

[58] E. F. Toro, M. Spruce, W. Speares, Restoration of the contact surface in the HLL-Riemann solver, Shock Waves 4 (1994) 25–34. doi:10.1007/BF01414629.

[59] H. Ranocha, M. Schlottke-Lakemper, A. R. Winters, E. Faulhaber, J. Chan, G. Gassner, Adaptive numerical simulations with Trixi.jl: A case study of Julia for scientific computing, Proceedings of the JuliaCon Conferences 1 (2022) 77. doi:10.21105/jcon.00077.

[60] M. Schlottke-Lakemper, A. R. Winters, H. Ranocha, G. Gassner, A purely hyperbolic discontinuous Galerkin approach for self-gravitating gas dynamics, Journal of Computational Physics 442 (2021) 110467. doi:10.1016/j.jcp.2021.110467.

[61] M. Schlottke-Lakemper, G. Gassner, H. Ranocha, A. R. Winters, J. Chan, Trixi.jl: Adaptive high-order numerical simulations of hyperbolic PDEs in Julia, https://github.com/trixi-framework/Trixi.jl, 2021. doi:10.5281/zenodo.3996439.

[62] P. J. van Der Houwen, B. P. Sommeijer, On the internal stability of explicit, m-stage Runge-Kutta methods for large m-values, ZAMM-Journal of Applied Mathematics and Mechanics/Zeitschrift für Angewandte Mathematik und Mechanik 60 (1980) 479–485. doi:10.1002/zamm.19800601005.

[63] J. G. Verwer, W. H. Hundsdorfer, B. P. Sommeijer, Convergence properties of the Runge-Kutta-Chebyshev method, Numerische Mathematik 57 (1990) 157–178. doi:10.1007/BF01386405.

[64] S. O'Sullivan, Factorized Runge–Kutta–Chebyshev methods, in: Journal of Physics: Conference Series, volume 837, IOP Publishing, 2017, p. 012020. doi:10.1088/1742-6596/837/1/012020.

[65] D. Doehring, G. J. Gassner, M. Torrilhon, Many-stage optimal stabilized Runge-Kutta methods for hyperbolic partial differential equations, Journal of Scientific Computing 99 (2024). doi:10.1007/s10915-024-02478-5.

[66] W. H. Hundsdorfer, J. G. Verwer, Numerical solution of time-dependent advection-diffusion-reaction equations, volume 33 of *Springer Series in Computational Mathematics*, 1 ed., Springer, 2010. doi:10.1007/978-3-662-09017-6.

[67] D. I. Ketcheson, Highly efficient strong stability-preserving Runge-Kutta methods with low-storage implementations, SIAM Journal on Scientific Computing 30 (2008) 2113–2136. doi:10.1137/07070485X.

[68] C. Rackauckas, Q. Nie, DifferentialEquations.jl – a performant and feature-rich ecosystem for solving differential equations in julia, The Journal of Open Research Software 5 (2017) 10. doi:10.5334/jors.151.

[69] A. J. Lotka, Undamped oscillations derived from the law of mass action., Journal of the American Chemical Society 42 (1920) 1595–1599. doi:10.1021/ja01453a010.

[70] V. Volterra, Variations and fluctuations of the number of individuals in animal species living together, ICES Journal of Marine Science 3 (1928) 3–51. doi:10.1093/icesjms/3.1.3, translated by M. E. Wells.

[71] J.-L. Boulnois, An exact closed-form solution of the Lotka-Volterra equations, 2023. doi:10.48550/arXiv.2303.09317.

[72] C. Burstedde, L. C. Wilcox, O. Ghattas, p4est: Scalable algorithms for parallel adaptive mesh refinement on forests of octrees, SIAM Journal on Scientific Computing 33 (2011) 1103–1133. doi:10.1137/100791634.

[73] D. A. Kopriva, A conservative staggered-grid Chebyshev multidomain method for compressible flows. II. A semi-structured method, Journal of Computational Physics 128 (1996) 475–488. doi:10.1006/jcph.1996.0225.

[74] S. F. Davis, Simplified second-order Godunov-type methods, SIAM Journal on Scientific and Statistical Computing 9 (1988) 445–473. doi:10.1137/0909030.

[75] E. F. Toro, Riemann solvers and numerical methods for fluid dynamics: A practical introduction, 3 ed., Springer Berlin, Heidelberg, 2013. doi:10.1007/b79761.

[76] M. H. Carpenter, C. A. Kennedy, Fourth-order 2N-storage Runge-Kutta schemes, Technical Report, NASA Langley Research Center, 1994. URL: https://ntrs.nasa.gov/api/citations/19940028444/downloads/19940028444.pdf.

[77] J. Niegemann, R. Diehl, K. Busch, Efficient low-storage Runge-Kutta schemes with optimized stability regions, Journal of Computational Physics 231 (2012) 364–372. doi:10.1016/j.jcp.2011.09.003.

[78] M. Parsani, D. I. Ketcheson, W. Deconinck, Optimized explicit Runge–Kutta schemes for the spectral difference method applied to wave propagation problems, SIAM Journal on Scientific Computing 35 (2013) A957–A986. doi:10.1137/120885899.

[79] T. Toulorge, W. Desmet, Optimal Runge–Kutta schemes for discontinuous Galerkin space discretizations applied to wave propagation problems, Journal of Computational Physics 231 (2012) 2067–2091. doi:10.1016/j.jcp.2011.11.024.

[80] H. Ranocha, L. Dalcin, M. Parsani, D. I. Ketcheson, Optimized Runge-Kutta methods with automatic step size control for compressible computational fluid dynamics, Communications on Applied Mathematics and Computation 4 (2022) 1191–1228. doi:10.1007/s42967-021-00159-w.

[81] T. Chen, C.-W. Shu, Entropy stable high order discontinuous Galerkin methods with suitable quadrature rules for hyperbolic conservation laws, Journal of Computational Physics 345 (2017) 427–461. doi:10.1016/j.jcp.2017.05.025.

[82] P. Chandrashekar, Kinetic energy preserving and entropy stable finite volume schemes for compressible Euler and Navier-Stokes equations, Communications in Computational Physics 14 (2013) 1252–1286. doi:10.4208/cicp.170712.010313a.

[83] F. Bassi, S. Rebay, A high-order accurate discontinuous finite element method for the numerical solution of the compressible Navier–Stokes equations, Journal of Computational Physics 131 (1997) 267–279. doi:10.1006/jcph.1996.5572.

[84] G. J. Gassner, A. R. Winters, F. J. Hindenlang, D. A. Kopriva, The BR1 scheme is stable for the compressible Navier–Stokes equations, Journal of Scientific Computing 77 (2018) 154–200. doi:10.1007/s10915-018-0702-1.

[85] C. A. Kennedy, M. H. Carpenter, R. M. Lewis, Low-storage, explicit Runge–Kutta schemes for the compressible Navier–Stokes equations, Applied numerical mathematics 35 (2000) 177–219. doi:10.1016/S0168-9274(99)00141-5.

[86] A. Uranga, P.-O. Persson, M. Drela, J. Peraire, Implicit large eddy simulation of transition to turbulence at low Reynolds numbers using a Discontinuous Galerkin method, International Journal for Numerical Methods in Engineering 87 (2011) 232–261. doi:10.1002/nme.3036.

[87] M. R. López-Morales, J. R. Bull, J. Crabill, T. D. Economon, D. E. Manosalvas, J. Romero, A. Sheshadri, J. E. Watkins, D. M. Williams, F. Palacios, A. Jameson, Verification and validation of HiFiLES: a high-order LES unstructured solver on multi-GPU platforms, in: 32nd AIAA Applied Aerodynamics Conference, 2014, pp. 1–27. doi:10.2514/6.2014-3168.

[88] A. Powell, Theory of Vortex Sound, The Journal of the Acoustical Society of America 36 (1964) 177. doi:10.1121/1.1918931.

[89] D. J. Lee, S. O. Koo, Numerical study of sound generation due to a spinning vortex pair, AIAA Journal 33 (1995) 20–26. doi:10.2514/3.12327.

[90] Michael Schlottke-Lakemper, Hans Yu, Sven Berger, Matthias Meinke, Wolfgang Schröder, A fully coupled hybrid computational aeroacoustics method on hierarchical Cartesian meshes, Computers & Fluids 144 (2017) 137–153. doi:10.1016/j.compfluid.2016.12.001.

[91] M. Schlottke-Lakemper, A direct-hybrid method for aeroacoustic analysis, Ph.D. thesis, RWTH Aachen University, 2017. URL: https://publications.rwth-aachen.de/record/688887/files/688887.pdf. doi:10.18154/RWTH-2017-04082.

[92] F. Alizadeh, D. Goldfarb, Second-order cone programming, Mathematical Programming, Series B 95 (2002) 3–51. doi:10.1007/s10107-002-0339-5.

[93] S. Boyd, L. Vandenberghe, Convex optimization, Cambridge University Press, 2004. URL: https://web.stanford.edu/~boyd/cvxbook/bv_cvxbook.pdf.

[94] S. Niederberger, A. Domahidi, E. Chu, S. Boyd, EiCOS: Eigen Conic Solver, https://github.com/EmbersArc/EiCOS, 2022. Accessed: 2022-06-16. Forked from commit: 87398a9c4732b09c5b2524470b67ad583719423b.

[95] A. Domahidi, E. Chu, S. Boyd, ECOS: An SOCP solver for embedded systems, in: European Control Conference (ECC), 2013, pp. 3071–3076.

[96] M. Grant, S. Boyd, CVX: Matlab software for disciplined convex programming, version 2.1, https://cvxr.com/cvx, 2014.

[97] M. Grant, S. Boyd, Graph implementations for nonsmooth convex programs, in: V. Blondel, S. Boyd, H. Kimura (Eds.), Recent Advances in Learning and Control, Lecture Notes in Control and Information Sciences, Springer-Verlag Limited, 2008, pp. 95–110. http://stanford.edu/~boyd/graph_dcp.html.

[98] G. Guennebaud, B. Jacob, et al., Eigen v3, http://eigen.tuxfamily.org, 2010.

[99] Boost, Boost C++ Libraries, Version 1.74.0.3ubuntu7, http://www.boost.org/, 2024. Last accessed 2024-02-21.

[100] D. Doehring, SSOCs: Arbitrary precision stability polynomials by second-order cones, https://github.com/DanielDoehring/SSOCs, 2024. doi:10.5281/zenodo.11184359, GitHub repository.